\documentclass[a4paper, 12pt]{elsarticle}
\usepackage{enumerate}
\usepackage{lineno,hyperref}
\usepackage{lineno,hyperref,tikz,adjustbox}
\usepackage{tikz-cd}
\tikzcdset{scale cd/.style={every label/.append style={scale=#1},
		cells={nodes={scale=#1}}}}
	
\modulolinenumbers[0]

\journal{Journal of \LaTeX\ Templates}

\usepackage{amssymb}
\usepackage{amsmath}
\usepackage{multirow}
\usepackage{longtable}
\usepackage{comment}
\usepackage{placeins}
\usepackage{mathtools}
\usepackage{algorithm}
\usepackage{algpseudocode}
\usepackage{enumitem}
\usepackage[utf8]{inputenc}
\usepackage{booktabs}
\usepackage{array}

\usepackage[draft,nomargin,inline]{fixme}
\fxsetface{inline}{\itshape}
\fxsetface{env}{\itshape}
\fxusetheme{color}

\usepackage{url}
\urldef{\mailsa}\path|{djukanovic, raidl}@ac.tuwien.ac.at,|
\urldef{\mailsb}\path|christian.blum@iiia.csic.es|

\usepackage{tikz}
\usetikzlibrary{positioning}
\definecolor{canaryyellow}{rgb}{1.0, 0.94, 0.0}
\definecolor{brightgreen}{rgb}{0.4, 1.0, 0.0}
\definecolor{jazzberryjam}{rgb}{0.65, 0.04, 0.37}

\usepackage{bookmark}
\usepackage{hypcap} 
\evensidemargin\oddsidemargin
\usepackage{graphicx}
\usepackage{xcolor}

\usepackage{xspace}
\usepackage{color}
\usepackage{epsfig}
\usepackage{caption}
\usepackage{subcaption}
\usepackage{mathrsfs}
\usepackage{amssymb}
\usepackage{amsmath}
\usepackage{amsthm}

\usepackage{blindtext}

\newcommand\undermat[2]{%
  \makebox[0pt][l]{$\smash{\underbrace{\phantom{%
    \begin{matrix}#2\end{matrix}}}_{\text{$#1$}}}$}#2}

\algnewcommand\algorithmicforeach{\textbf{for each}}
\algdef{S}[FOR]{ForEach}[1]{\algorithmicforeach\ #1\ \algorithmicdo}

\newtheorem{Theorem}{Theorem}
\newtheorem{Lemma}{Lemma}
\newtheorem{Proposition}{Proposition}
\newtheorem{Observation}{Observation}
\newtheorem{Corollary}{Corollary}
\newtheorem{Definition}{Definition}

\begin{document}

\begin{frontmatter}

\title{Theoretical Studies of the $k$-Strong Roman Domination Problem}

\author[1]{Bojan Nikoli\'c\fnref{myfootnote}}
\address[1]{  Faculty of Natural Sciences and Mathematics, University of Banja Luka, Mladena Stojanovi\'ca 2, 78000 Banja Luka, Bosnia and Herzegovina}
\ead[2]{bojan.nikolic@pmf.unibl.org}
\author[1]{Marko Djukanovi\'c}
\fntext[myfootnote]{Corresponding author}
\ead[2]{marko.djukanovic@pmf.unibl.org}
\author[1]{Milana Grbi\'c}
\ead[3]{milana.grbic@pmf.unibl.org}
\author[1]{Dragan Mati\'c}
\ead[4]{dragan.matic@pmf.unibl.org}

\begin{abstract}
     The concept of Roman domination has been a subject of intrigue for more than two decades with the fundamental Roman domination problem standing out as one of the most significant challenges in this field. This article studies a practically motivated generalization of this problem, known as the \emph{$k$-strong Roman domination}. In this variation, defenders within a network are tasked with safeguarding any $k$ vertices simultaneously,  under multiple attacks.  The objective is to find a feasible mapping that assigns
      an (integer) weight to each vertex of the input graph with a minimum sum of weights across all vertices. A function is considered feasible if any non-defended vertex, i.e. one labeled by zero, is protected by at least one of its  neighboring vertices labeled by at least two. Furthermore, each defender ensures the safety of a non-defended vertex by imparting a value of one to it while always retaining a one for themselves.  To the best of our knowledge, this paper represents  the first theoretical study on this problem. The study presents results for general graphs, establishes connections between the problem at hand and other domination problems, and provides exact values and bounds for specific graph classes, including complete graphs, paths, cycles, complete bipartite graphs, grids, and a few selected classes of convex polytopes. Additionally, an attainable lower bound for general cubic graphs is provided.
\end{abstract}


\end{frontmatter}


 	\section{Introduction}\label{sec:introduction}

     Within the last few decades, the dedication to solving various domination problems on graphs theoretically and also computationally was highly intensified, see, for example, surveys~\cite{cockayne1978domination,henning2009survey, goddard2013independent,lu2010survey,chang1998algorithmic}.
   Many variants of the basic domination problem on graphs have found a wide range of practical applications in telecommunications~\cite{balasundaram2006graph},  scheduling problems~\cite{rao2021survey}, detecting central proteins in biological networks~\cite{milenkovic2011dominating}, having also a proven relation to facility location problems~\cite{gupta2013domination}, and many others.
   Among these variants, the \emph{Roman domination problem} (RDP) has been placed high in the literature as it has intrigued many researchers to investigate it from both perspectives, theoretically and practically, see~\cite{chellali2020roman,cockayne2004roman}. More in detail, this problem concerns an optimal assigning of the values from $\{0, 1, 2\}$  to the vertices of a graph fulfilling the constraint that each \emph{weak} vertex, i.e. one labelled with 0,  has to have at least one neighboring vertex assigned with 2, capable to protect it in the case of an attack. The price of such an assigning on the graph is given by summing up all labels of the vertices of that graph.  The practical application of RDP
   has been in discovering an optimal military strategy~\cite{henning2003defending}, network design of server placements~\cite{pagourtzis2002server}, etc. The RDP resembles the well-known maximal covering location problem~\cite{church1974maximal}.  Concerning the theoretical studies of RDP, the exact bounds are found and proven for various graph classes such as paths, complete graphs, grid graphs, complete bipartite graphs, convex polytopes,  etc.~\cite{cockayne2004roman,kartelj2021roman, liu2013roman}.  For having a look at some variants of RDP so far studied in the literature, we refer readers   to~\cite{chellali2021varieties}.
   Among dozens of variants of RDP which appear in literature, we mention the \emph{Roman $k$-domination problem}~\cite{kammerling2009roman}, related to the problem we study here. A Roman $k$-dominating function on graph $G$ is a function $f \colon V (G) \rightarrow \{0, 1, 2\}$
   such that every vertex labelled with 0 is adjacent to at least $k$ vertices labelled with 2.
   The minimum weight of a Roman $k$-dominating function on a graph $G$ is called the Roman $k$-domination number denoted as $\gamma_{kR}(G)$.
 	
 	In this work, we study another practically motivated variant of RDP, so-called the  \emph{k-Strong Roman Domination} ($k$-SRD) \emph{Problem}~\cite{liu2020k}. It represents a generalization of the basic RDP, where instead of having attacked one vertex at the same time, here any $k\geqslant 1$ vertices are attacked simultaneously and all such attack possibilities need to be defended by itself or some neighboring vertices. We again ask for an optimal assignment fulfilling this demand regarding all attacks. 	Note that the 1-SRD problem gives the basic RDP.  In reality, defenders often have to be prepared for multiple simultaneous attacks to ensure defending more than one location at the same time. One such example is planning anti-terrorism strategies. Another possible application of $k$-SRD problem may be found in disaster relief, and supply chain disruptions~\cite{liu2020k}.   $k$-SRD problem is indeed theoretically much harder to solve than so far considered variants of RD problem as it belongs to the complexity class EXPTIME~\cite{schoning1986complete}.  Thus, being quite challenging, it is expected to be computationally (practically as well) much harder to solve than most of the other RD-like problems known in the literature.  As we will see, the proofs dedicated to the $k$-SRD problem appear to be usually more complex than the corresponding ones for other RD-like problems considered on the same graph classes. It is worth mentioning a problem related to $k$-SRD, the \emph{strong Roman domination} problem~\cite{alvarez2017strong},  that also considers multiple simultaneous attacks on graphs under different constraints. In this variant of the Roman domination problem, the vertex whose label is greater than 1 is considered to be ``strong'' and besides itself is also able to defend half of its neighbors. The $k$-SRD problem is introduced as relaxation of strong Roman domination problem~\cite{liu2020k}.

 Concerning the literature review about the $k$-SRD problem, we are aware of just one paper published so far~\cite{liu2020k}, where the problem is studied from the practical perspective; two algorithms are developed: an integer linear programming and a Benders-based decomposition method. Both of these methods are limited to the application on small-sized random graphs (up to 50 vertices and up to $k=5$). This work is the first of the kind in the literature which studies the $k$-SRD problem from the theoretical point of view.

    \subsection{Notation}
 	
 	Let $G=(V(G),E(G))$ be a simple, undirected graph where $V(G)$ and $E(G)$ represents its vertices and edges.
Let $|V(G)|$ be the number of vertices of graph $G$. By $d(v), v \in V(G)$, we denote the degree of vertex $v$, and by $\Delta(G)= \max  \{ d(v) \mid v \in V(G) \}$ the maximum degree of graph $G$.  With $N(v)$ we denote the set of adjacent vertices of vertex $v$.  Let $f$ be an (integer-domain) vertex mapping on $V(G)$. By $V_j$ we define the set of all vertices that are mapped to integer $j$, i.e. $V_j = \{ v \in V(G) \mid f(v) = j\}$, additionally, $V_{\geqslant j} = \{ v \in V(G) \mid f(v) \geqslant j\}, j \in \mathbb{N}$.
Notice that there is a 1-1 correspondence between $f:V\rightarrow \{0,1,\ldots,q\}$ and the ordered partition $(V_0,V_1,\ldots,V_q)$, so each  function $f$  can also be denoted as $f=(V_0,V_1,\ldots,V_q)$.
Last but not least,  $\omega(f, G) = \sum_{v \in V(G) } f(v)$. If not stated differently, $\omega(f, G)$ will be simply denoted by $\omega(f)$. We always consider $|V(G)|\geqslant  k$ throughout the paper.
 	
 	\vspace{0.5cm}

 	\subsection{Problem Definition}

 	For a graph $G$, let $P := \{v_1, v_2, \ldots, v_k\}$, $v_1, v_2, \ldots , v_k \in V(G)$,
 be an attack pattern, that is,  the set of vertices under attack. In essence,  for given graph $G$ and given value $k$ there are $\binom{|V(G)|}{k}$ possible attack patterns. By $\mathcal{P}(G)= \{ P_1, \ldots, P_{\binom{|V(G)|}{k}}\}$ we denote the set of the enumerated attack patterns in graph $G$.

 \Definition\label{def:SRDfunction} A function $f: V(G) \mapsto \{0,1, \ldots, \min(\Delta(G), k) + 1\}$ is a proper $k$-SRD function if, for each attack pattern $P \in \mathcal{P}(G)$, the following is fulfilled:
 \begin{itemize}
 	\item [(i)] Every vertex labeled with $f(v)\geqslant 2$ can defend at most $f(v)-1$   neighboring vertices from $P$ labeled with 0;
 	\item[(ii)] Vertices labelled with 0 belonging to the attack pattern $P$ have to be defended by at least one of its neighboring vertices.
 \end{itemize}
 \normalfont

 Note that if all vertices are labelled with 1, then all of them are ``safe'' under any attack pattern from $\mathcal{P}(G)$, thus obtaining a (trivial) $k$-SRD function on a graph $G$.

 The problem of $k$-SRD  on graph $G$ is to find a proper $k$-SRD function of a minimum cumulative weight, called $k$-SRD number, denoted by $\gamma_{k-\text{SRD}}(G)$, i.e.  $\gamma_{k-\text{SRD}}(G) := \min  \{\omega(f) \mid f \text{ is a proper } k\text{-SRD on graph } G  \}$.  Such a $k$-SRD function for which $\gamma_{k-\text{SRD}}(G)$ is reached we call $\gamma_{k-\text{SRD}}$ function.

 Since any vertex $v\in V(G)$ cannot protect more vertices than $d(v)$+1, the maximum value of any $\gamma_{k-\text{SRD}}$  function is clearly not greater than $\min(\Delta(G), k) + 1$ \cite{liu2020k}.

 	Let $f \colon V(G) \mapsto \{0,1, \ldots, \min(\Delta(G), k) + 1\}$ be a candidate  $k$-SRD function. An attack pattern containing only vertices labeled with a positive value are considered protected, as each such vertex is self-protected. Therefore, we are interested only in attack patterns that contain at least one vertex labeled with zero.
 	On the other hand, any subset $B\subset V_0$, $0<|B|\leqslant k$ belongs to at least one attack pattern (more precisely, it belongs to exactly $\binom{k}{|B|}$ attack patterns). For each such set $B$, there exists a set $A\subset V_{\geqslant 2}$, such that each vertex from $B$ is protected by a vertex from $A$, following the conditions of Definition~\ref{def:SRDfunction}. Therefore, it implies $\sum_{v\in A}f(v)  \geqslant |A|+|B|$.

 \subsection{Contributions}
 	
  Contributions to the paper are as follows:
\begin{itemize}
	
	\item We present a few results for the case of fixed $k=2$, proving the attainable lower bound for general cubic graphs using  the well-known discharging approach~\cite{shao2018discharging}. Furthermore, exact bounds are derived for specific cubic graphs, specifically targeting a few classes of convex polytopes.
	
	\item  For general graphs, relations on $k$-SRD number and the optimal values of some related problems such as the Roman $k$-domination and the basic domination problem are stated and formally proved. Some basic properties of $k$-SRD number as a function depending on the value of $k$ or general input graph $G$ are proved.
	
   \item Exact values of $k$-SRD numbers are determined and proven for various graph classes, including complete graphs, paths, star graphs, wheels, complete bipartite graphs, and grid graphs of size $2 \times n$, all for arbitrary $k\geqslant 2$. In the proof of lower bounds of some other graph classes (like path and grid graphs), we propose a new memory-based
       discharging procedure. This novel technique improves the ``basic" discharging procedure in a way that common defenders of the same zero-labeled vertex share the information
       about the amount of charge that each one of them is giving to this vertex. Essentially, this approach enables more flexibility in weights redistribution when comes to a
       multi-member attack pattern scenario. 

\end{itemize}
 	
 	 \section{Main Results}

      Section~\ref{sec:general-graph-ksrd-k-2} provides the results about $\gamma_{2-SRD}$ number on cubic graphs and a few special graphs classes  of convex polytopes. Section~\ref{sec:general-graph-ksrd} presents the results about the exact, lower and upper bounds on  $\gamma_{k-\text{SRD}}$ number for general $k \geqslant 2$,  for several graph classes including path graphs, complete graphs, complete bipartite graphs and grid graphs.

\subsection{Results for $2$-SRD problem}\label{sec:general-graph-ksrd-k-2}

In this section we provide some bounds for the special case  2-SRD problem. We start with a lower bound for cubic graphs whose all vertices are of degree three.

Let $G=(V(G),E(G))$ be a cubic graph and  let $f$ be an  arbitrary ${2}$-SRD function. Since $\Delta(G)=3$, one can consider that  $f = (V_{0},V_1,V_2,V_3)$.
The following observation holds.
\Observation \label{obs:2}  A vertex from $V_2$ can have at most one adjacent vertex from $V_0$ which is not adjacent with another vertex from $V_2\cup V_3$. In other words, for a given attack pattern, two vertices from $V_0$ can not be solely protected by the same vertex from $V_2$.

\Theorem \label{thm:cubic-lower-bound} Let $G=(V(G),E(G))$ be a cubic graph and $|V(G)|=n$. Then

$$\gamma_{2-SRD}(G)\geqslant \frac{2n}3$$

\proof
Let $f$ be an arbitrary $2$-SRD function and $f=(V_0,V_1,V_2,V_3)$.

We define  function $g$ by the following  basic discharging procedure  rules.
\begin{enumerate}
	\item [{\bf R0}] For each $v\in V(G)$, we set $g(v) = f(v)$ at the beginning of the procedure.
	\item[{\bf R1}]  If $u\in V_0$ and $v\in N(u)\cap V_3$, then vertex $u$ gets charge of $\frac23$ from vertex $v$, i.e. we increase the value $g(u)$ by $\frac{2}{3}$
     and also decrease the value $g(v)$ by $\frac{2}{3}$.
	\item [{\bf R2}] If $u\in V_0$ and $N(u)\cap V_2=\{v\}$, then vertex $u$ gets charge of $\frac23$ from vertex $v$, i.e. we increase the value $g(u)$ by $\frac{2}{3}$
     and also decrease the value $g(v)$ by $\frac{2}{3}$.
	\item [{\bf R3}] If $u\in V_0$ and $v,z\in N(u)\cap V_2$, then vertex $u$ gets charge of $\frac13$ from each of vertices $v,z$, i.e. we increase the value $g(u)$ by $\frac{2}{3}$ and also decrease each of the values $g(v)$ and $g(z)$ by $\frac{1}{3}$.
\end{enumerate}

Notice that each rule of this procedure preserves the cumulative weight of the function $f$, implying $\omega(g) = \omega(f)$.
We claim that after the application of the previous rules, for any vertex $v\in V(G)$, it holds $g(v)\geqslant \frac{2}{3}$. Indeed,

\begin{itemize}
	\item For each $v\in V_1$, $g(v) = f(v) = 1$.
	\item Let $v\in V_3$. Since $v$ can have at most 3 adjacent vertices from $V_0$, then  $g(v) \geqslant 3-3\cdot \frac{2}3 \geqslant 1$.
	
	\item Let $v\in V_0$. According to Rules R1-R3, obviously $g(v)\geqslant \frac23$.
	
	\item Let $v\in V_2$.
	According to Observation~\ref{obs:2} and since $d(v)=3$, $v$ can  charge $\frac23$ to a vertex from $V_0$ at most once and can charge $\frac13$ at most twice.  Therefore,  $g(v)\geqslant 2-\frac23 - 2\cdot \frac 13 = \frac 23$.
\end{itemize}

Therefore, $\omega(f) = \omega(g) \geqslant \frac 2 3 n$, which concludes the proof.
\qed \vspace{0.2cm}

The next three propositions are dedicated to the cubic graph classes whose   (exact) $2$-SRD numbers attain the theoretical lower bound derived in the last theorem.

The convex polytope $D_n,n\geqslant 5$~\cite{filipovic2022solving}, consists of $2n$ 5-sided faces and 2 $n$-sided faces.  Formally, graph $D_n=(V(D_n), E(D_n))$ consists of the following sets
of vertices and edges:
\begin{align*}
	&V(D_n) = \{ a_i, b_i, c_i, d_i \mid i \in \{0, \ldots, n-1\}\}, \\
	&V(D_n) = \{ a_ia_{i+1}, d_id_{i+1},a_ib_i, b_ic_i, b_{i+1}c_i, c_id_i  \mid i \in \{0, \ldots, n-1\}  \}.
\end{align*}

\begin{figure}[H]
	\centering
	\setlength\unitlength{1mm}
	\begin{picture}(130,60)
		\thicklines
		\tiny
		\put(47.7,5.0){\circle*{1}} \put(47.7,5.0){\line(1,1){7.3}} \put(47.7,5.0){\line(-5,3){8.7}}
		\put(55.0,12.3){\circle*{1}} \put(55.0,12.3){\line(4,1){10.0}} \put(55.0,12.3){\line(-3,5){5.0}}
		\put(65.0,15.0){\circle*{1}} \put(65.0,15.0){\line(4,-1){10.0}} \put(65.0,15.0){\line(0,1){10.0}}
		\put(75.0,12.3){\circle*{1}} \put(75.0,12.3){\line(1,-1){7.3}} \put(75.0,12.3){\line(3,5){5.0}}
		\put(82.3,5.0){\circle*{1}} \put(82.3,5.0){\line(5,3){8.7}}
		\put(39.0,10.0){\circle*{1}} \put(39.0,10.0){\line(-5,-2){12.7}} \put(39.0,10.0){\line(-1,6){2.3}}
		\put(50.0,21.0){\circle*{1}} \put(50.0,21.0){\line(-6,1){13.3}} \put(50.0,21.0){\line(2,5){4.6}}
		\put(65.0,25.0){\circle*{1}} \put(65.0,25.0){\line(-5,4){10.4}} \put(65.0,25.0){\line(5,4){10.4}}
		\put(80.0,21.0){\circle*{1}} \put(80.0,21.0){\line(-2,5){4.6}} \put(80.0,21.0){\line(6,1){13.3}}
		\put(91.0,10.0){\circle*{1}} \put(91.0,10.0){\line(1,6){2.3}} \put(91.0,10.0){\line(5,-2){12.7}}
		\put(26.4,5.4){\circle*{1}} \put(26.4,5.4){\line(-4,1){9.7}}
		\put(36.7,23.3){\circle*{1}} \put(36.7,23.3){\line(-1,1){7.1}}
		\put(54.6,33.6){\circle*{1}} \put(54.6,33.6){\line(-1,4){2.6}}
		\put(75.4,33.6){\circle*{1}} \put(75.4,33.6){\line(1,4){2.6}}
		\put(93.3,23.3){\circle*{1}} \put(93.3,23.3){\line(1,1){7.1}}
		\put(103.6,5.4){\circle*{1}} \put(103.6,5.4){\line(4,1){9.7}}
		\put(16.7,7.9){\circle*{1}} \put(16.7,7.9){\line(3,5){12.9}}
		\put(29.6,30.4){\circle*{1}} \put(29.6,30.4){\line(5,3){22.4}}
		\put(52.1,43.3){\circle*{1}} \put(52.1,43.3){\line(1,0){25.9}}
		\put(77.9,43.3){\circle*{1}} \put(77.9,43.3){\line(5,-3){22.4}}
		\put(100.4,30.4){\circle*{1}} \put(100.4,30.4){\line(3,-5){12.9}}
		\put(113.3,7.9){\circle*{1}}
		\thinlines
		\put(47.7,5.0){\line(-1,-2){3}} \put(82.3,5.0){\line(1,-2){3}}
		\put(26.4,5.4){\line(1,-1){4}} \put(103.6,5.4){\line(-1,-1){4}}
		\put(16.7,7.9){\line(0,-1){6}} \put(113.3,7.9){\line(0,-1){6}}
		
		\put(48.3,3.5){$a_2$} \put(54.5,9.7){$a_1$} \put(63.0,12.0){$a_0$} \put(71.5,9.7){$a_{n-1}$} \put(77.7,3.5){$a_{n-2}$}
		\put(36.3,11.0){$b_2$} \put(48.0,22.7){$b_1$} \put(64.0,27.0){$b_0$} \put(80.0,22.7){$b_{n-1}$} \put(91.7,11.0){$b_{n-2}$}
		\put(24.4,8.7){$c_2$} \put(36.4,26.2){$c_1$} \put(55.5,35.3){$c_0$} \put(76.7,33.6){$c_{n-1}$} \put(94.2,21.6){$c_{n-2}$} \put(103.3,2.5){$c_{n-3}$}
		\put(13.8,8.5){$d_2$} \put(27.2,31.8){$d_1$} \put(50.5,45.2){$d_0$} \put(77.5,45.2){$d_{n-1}$} \put(100.8,31.8){$d_{n-2}$} \put(114.2,8.5){$d_{n-3}$}
		
	\end{picture}
	\caption{The graph of convex polytope $D_n$.} \label{fig:d_n}
\end{figure}
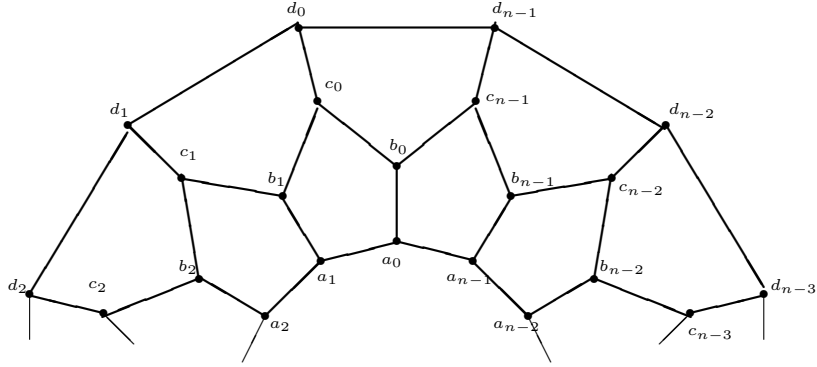

Note that indices of vertices are using the notation by modulo $n$; for instance $a_n = a_0$. The graph of convex polytope $D_n, n \geqslant 5$, is given in Figure~\ref{fig:d_n}.

\Proposition For $n\geqslant 5$, it holds
$$\gamma_{2-SRD}(D_n) =  \left\lceil \frac{2|V(D_n)|}{3}\right\rceil =  \left\lceil \frac{8n}{3} \right\rceil. $$


\proof

We define function $f:V(D_n)\mapsto\{0,1,2\}$ according to the following cases:

\begin{itemize}

\item $n \equiv 0 (\bmod\ 3)$;

For $1\leqslant i\leqslant n-1$, we put  $$ f(a_i) = f(c_i)  = \begin{cases}
	2, i \equiv 1 (\bmod\ 3) \\
	0, \textrm{otherwise}
\end{cases},
f(b_i) = f(d_i)=\begin{cases}
	2,       i \equiv 0(\bmod\ 3)\\
	0, \textrm{otherwise}.
\end{cases}$$

One can easily see that
$\omega(f)=  2 \cdot \frac{n}{3} +  2 \cdot \frac{n}{3} +  2 \cdot \frac{n}{3} +  2 \cdot \frac{n}{3} = \frac{8n}{3}=\lceil \frac{8n}{3}\rceil$.


\item $n \equiv 1 (\bmod\ {3})$;

For $0\leqslant i< n-1$, we put
$$ f(a_i) = f(c_i)  = \begin{cases}
	2, i \equiv 1 (\bmod\ 3) \\
	0, \textrm{otherwise}
\end{cases},
f(b_i)=\begin{cases}
	2,       i \equiv 0(\bmod\ 3)\\
	0, \textrm{otherwise}.
\end{cases}$$
For $1\leqslant i< n-1$, we put
$$
f(d_i) = \begin{cases}
	2, i \equiv 2 (\bmod\ 3) \\
	0,  \textrm{otherwise.}
\end{cases}
$$

The remaining vertices are labeled as follows:
$$f(a_{n-1})=f(c_{n-1}) = f(d_{n-1})=0,\ f(b_{n-1})= 2, f(d_0) = 1.$$
In this case, we have
$\omega(f)= 2\lfloor \frac{n}{3} \rfloor + 2\lfloor \frac{n}{3}\rfloor + 2 \lfloor \frac{n}{3} \rfloor + (2 \lfloor \frac{n}{3} \rfloor+1) = 8 \lfloor \frac{n}{3} \rfloor + 1 = \lceil \frac{8n}{3}\rceil$;

\item $n \equiv 2 (\bmod\ 3)$;

For $3 \leqslant i < n-2$, we put
$$f(a_i) =   \begin{cases}
	2, i \equiv 1(\bmod\ 3) \\
	0, \textrm{otherwise}
\end{cases},
f(b_i)= f(d_i) =  \begin{cases}
	2, i \equiv 2(\bmod\ 3) \\
	0, \textrm{otherwise}
\end{cases}$$
$$ f(c_i) = \begin{cases}
	2, i \equiv 0(\bmod\ 3) \\
	0, \textrm{otherwise}
\end{cases}
$$
The remaining vertices are labeled by
\begin{align*}
	&f(a_0)=0, f(a_1) = 2, f(a_2) = 0, f(a_{n-2})= 0, f(a_{n-1})= 2;\\
	&f(b_0) = 2, f(b_1)= f(b_2)  = f(b_{n-2}) = f(b_{n-1})= 0;\\
	&f(c_0) = 0, f(c_1) = f(c_2) =  f(c_{n-2})=2, f(c_{n-1})=0;\\
	&f(d_0)=2, f(d_1)= f(d_2)= f(d_{n-2})= f(d_{n-1})=0.
\end{align*}
In this case, we have $\omega(f)= (2\lfloor \frac{n}{3} \rfloor + 2) +2 \lfloor \frac{n}{3} \rfloor + (2\lfloor \frac{n}{3} \rfloor +4 ) + 2\lfloor \frac{n}{3} \rfloor = 8 \lfloor \frac{n}{3} \rfloor +6 = 8 \frac{n-2}{3} + 6=  \lceil \frac{8n}{3} \rceil$.

\end{itemize}

\begin{table}[H]
	\centering
	\begin{tabular}{c|l |c } \hline
		$v \in V_0$ & $N(v) \cap V_2$  & Indices  \\ \hline
		$a_{3i}$    & $a_{3i+1}$, $b_{3i}$    &  \\
		$a_{3i+2}$      &  $a_{3i+1}$   &  \\

		$c_{3i}$      &   $b_{3i}$, $d_{3i}$ & \\
		$c_{3i+2}$      &   $b_{3i+3}$  & $i=0, \ldots, \lfloor \frac{n}{3}\rfloor-1$ \\
		
		$b_{3i+1}$      &   $c_{3i+1}$, $a_{3i+1}$  &  \\
		$b_{3i+2}$      &  $c_{3i+1}$     &  \\

		$d_{3i+1}$      &   $d_{3i}$,  $c_{3i+1}$  &  \\
		$d_{3i+2}$      &  $d_{3i+3}$   &   \\

		\hline
	\end{tabular}
	\caption{Coverage of  $V_0$ vertices of function $f$ defined for $D_n$, the case $n \equiv 0$ ($\bmod\ 3$) } \label{opt-srdf-dn-3l}
	
\end{table}

\begin{table}[H]
	\centering
	\begin{tabular}{c|l |c } \hline
		$v \in V_0$ & $N(v) \cap V_2$  & Indices  \\ \hline
		$a_{3i}$    & $a_{3i+1}$, $b_{3i}$    & \\
		$a_{3i+2}$      &  $a_{3i+1}$   & \\

		$c_{3i}$      &   $b_{3i}$&   \\
		$c_{3i+2}$      &   $b_{3i+3}$, $d_{3i+2}$  & $i=0, \ldots, \lfloor \frac{n}{3}\rfloor-1$  \\
		
		$b_{3i+1}$      &   $c_{3i+1}$, $a_{3i+1}$  & \\
		$b_{3i+2}$      &  $c_{3i+1}$     &\\

		$d_{3i+1}$      &   $d_{3i+2}$,  $c_{3i+1}$  &\\\hline
		$d_{3i}$      &  $d_{3i-1}$   & $i=1, \ldots, \lfloor \frac{n}{3}\rfloor-1$  \\\hline
		$a_{n-1}$	&$b_{n-1}$&\\
		$c_{n-1}$	&$b_0,\ b_{n-1}$&\\\hline
		
		\hline
	\end{tabular}
	\caption{Coverage of  $V_0$ vertices of function $f$ defined for $D_n$, the case $n \equiv 1$ ($\bmod\ 3$) } \label{opt-srdf-dn-3l+1}
	
\end{table}

Tables ~\ref{opt-srdf-dn-3l}--\ref{opt-srdf-dn-3l+2} depict the coverage of $V_0$ vertices of the function $f$ with their neighbors from $V_2$. From these tables,
we can conclude the following:
\begin{itemize}
	\item Each vertex from $V_0$ is protected by at least one vertex from $V_2$;
	\item Each vertex from $V_2$ protects three vertices from $V_0$, but it protects only one vertex from $V_0$ solely.
\end{itemize}
Therefore,  for each case, any of $|V_0|\choose 2$ possible attack patterns can be defended, implying  $f$ is a 2-SRD function and $\gamma_{2-SRD}(D_n) \leqslant \lceil \frac{8n}{3} \rceil,$ for all $n \geqslant 5$. From this and Theorem~\ref{thm:cubic-lower-bound}, one can conclude that $\gamma_{2-SRD}(D_n)= \left \lceil  \frac{2}{3} |V(D_n)|\rceil = \right \lceil \frac{8n}{3} \rceil.$

\qed \vspace{0.2cm}

\begin{table}[H]
	\centering
	\begin{tabular}{c|l |c } \hline
		$v \in V_0$ & $N(v) \cap V_2$  & Indices  \\ \hline
		
		$a_{3i+2}$      &  $a_{3i+1}$, , $b_{3i+2}$   &  $i=0, \ldots, \lfloor \frac{n}{3}\rfloor$ \\\hline
		$a_{3i}$    & $a_{3i+1}$    &  \\
		
		$b_{3i+1}$      &   $a_{3i+1}$, $c_{3i}$  &  \\
		
		$c_{3i+1}$      &   $b_{3i+2}$& $i=1, \ldots, \lfloor \frac{n}{3}\rfloor$  \\
		$c_{3i+2}$      &   $b_{3i+2}$, $d_{3i+2}$  & \\
		$d_{3i+1}$      &   $d_{3i+2}$ &  \\\hline
		
		$b_{3i}$      &  $c_{3i}$     & $i=2, \ldots, \lfloor \frac{n}{3}\rfloor$  \\

		$d_{3i}$      &  $d_{3i-1}$, $c_{3i}$   &   \\\hline
		$a_0$&$a_1$, $a_{n-1}$, $b_0$&\\
		$b_1$&$a_1$, $c_1$&\\
		$b_2$&$c_1$, $c_2$&\\
		$b_3$&$c_2$, $c_3$&\\
		$c_0$&$b_0$, $d_0$&\\
		$d_1$&$c_1$, $d_0$&\\
		$d_2$&$c_2$&\\
		$d_3$&$c_3$&\\
		
		\hline
	\end{tabular}
	\caption{Coverage of  $V_0$ vertices of function $f$ defined for $D_n$, the case $n \equiv 2$ ($\bmod\ 3$) } \label{opt-srdf-dn-3l+2}
	
\end{table}

The graph of convex polytope $R_n''=(V(R_n''),E(R_n'')), n \geqslant 5$, introduced in~\cite{macdougall2006vertex}, consists of the
following sets of vertices and edges:
\begin{align*}
	&V(R_n'') = \{a_i, b_i, c_i, d_i, e_i, f_i \mid i \in \{0, \ldots, n-1 \}  \}, \\
	&E(R_n'') = \{ a_i a_{i+1}, f_i f_{i+1}, a_i b_i,  b_i c_i, c_id_i,
	d_ie_i, e_if_i, b_{i+1} c_i, d_{i} e_{i+1}  \mid i \in \{0, \ldots, n-1 \} \}.
\end{align*}

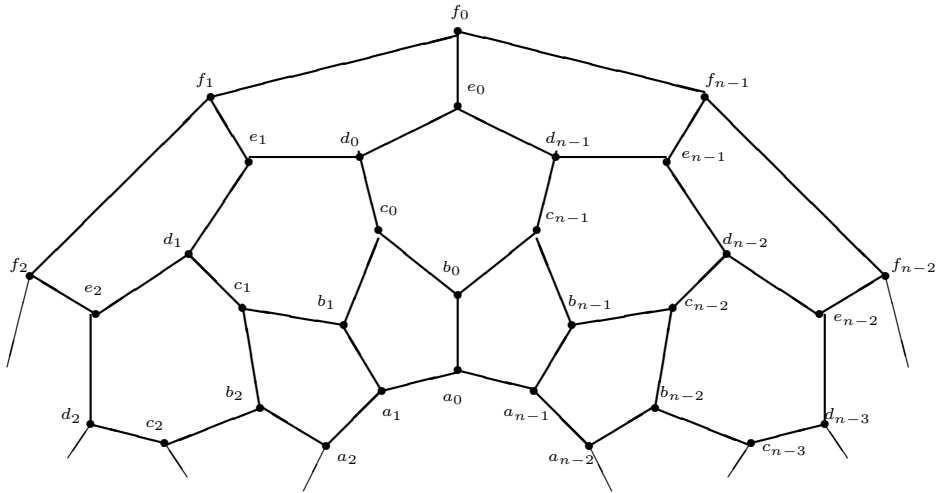
\begin{figure}[H]
	\centering

\setlength\unitlength{1mm}
\begin{picture}(130,60)
	\thicklines
	\tiny
	\put(47.7,5.0){\circle*{1}} \put(47.7,5.0){\line(1,1){7.3}} \put(47.7,5.0){\line(-5,3){8.7}}
	\put(55.0,12.3){\circle*{1}} \put(55.0,12.3){\line(4,1){10.0}} \put(55.0,12.3){\line(-3,5){5.0}}
	\put(65.0,15.0){\circle*{1}} \put(65.0,15.0){\line(4,-1){10.0}} \put(65.0,15.0){\line(0,1){10.0}}
	\put(75.0,12.3){\circle*{1}} \put(75.0,12.3){\line(1,-1){7.3}} \put(75.0,12.3){\line(3,5){5.0}}
	\put(82.3,5.0){\circle*{1}} \put(82.3,5.0){\line(5,3){8.7}}
	\put(39.0,10.0){\circle*{1}} \put(39.0,10.0){\line(-5,-2){12.7}} \put(39.0,10.0){\line(-1,6){2.3}}
	\put(50.0,21.0){\circle*{1}} \put(50.0,21.0){\line(-6,1){13.3}} \put(50.0,21.0){\line(2,5){4.6}}
	\put(65.0,25.0){\circle*{1}} \put(65.0,25.0){\line(-5,4){10.4}} \put(65.0,25.0){\line(5,4){10.4}}
	\put(80.0,21.0){\circle*{1}} \put(80.0,21.0){\line(-2,5){4.6}} \put(80.0,21.0){\line(6,1){13.3}}
	\put(91.0,10.0){\circle*{1}} \put(91.0,10.0){\line(1,6){2.3}} \put(91.0,10.0){\line(5,-2){12.7}}
	\put(26.4,5.4){\circle*{1}} \put(26.4,5.4){\line(-4,1){9.7}}
	\put(36.7,23.3){\circle*{1}} \put(36.7,23.3){\line(-1,1){7.1}}
	\put(54.6,33.6){\circle*{1}} \put(54.6,33.6){\line(-1,4){2.6}}
	\put(75.4,33.6){\circle*{1}} \put(75.4,33.6){\line(1,4){2.6}}
	\put(93.3,23.3){\circle*{1}} \put(93.3,23.3){\line(1,1){7.1}}
	\put(103.6,5.4){\circle*{1}} \put(103.6,5.4){\line(4,1){9.7}}
	\put(16.7,7.9){\circle*{1}} \put(16.7,7.9){\line(0,1){14.6}}
	\put(29.6,30.4){\circle*{1}} \put(29.6,30.4){\line(-3,-2){12.3}} \put(29.6,30.4){\line(2,3){7.9}}
	\put(52.1,43.3){\circle*{1}} \put(52.1,43.3){\line(-1,0){14.6}} \put(52.1,43.3){\line(2,1){12.9}}
	\put(77.9,43.3){\circle*{1}} \put(77.9,43.3){\line(-2,1){12.9}} \put(77.9,43.3){\line(1,0){14.6}}
	\put(100.4,30.4){\circle*{1}} \put(100.4,30.4){\line(-2,3){7.9}} \put(100.4,30.4){\line(3,-2){12.3}}
	\put(113.3,7.9){\circle*{1}} \put(113.3,7.9){\line(0,1){14.6}}
	\put(17.4,22.5){\circle*{1}} \put(17.4,22.5){\line(-5,3){8.7}}
	\put(37.5,42.6){\circle*{1}} \put(37.5,42.6){\line(-3,5){5.0}}
	\put(65.0,50.0){\circle*{1}} \put(65.0,50.0){\line(0,1){10.0}}
	\put(92.5,42.6){\circle*{1}} \put(92.5,42.6){\line(3,5){5.0}}
	\put(112.6,22.5){\circle*{1}} \put(112.6,22.5){\line(5,3){8.7}}
	\put(8.7,27.5){\circle*{1}} \put(8.7,27.5){\line(1,1){23.8}}
	\put(32.5,51.3){\circle*{1}} \put(32.5,51.3){\line(4,1){32.5}}
	\put(65.0,60.0){\circle*{1}} \put(65.0,60.0){\line(4,-1){32.5}}
	\put(97.5,51.3){\circle*{1}} \put(97.5,51.3){\line(1,-1){23.8}}
	\put(121.3,27.5){\circle*{1}}
	\thinlines
	\put(47.7,5.0){\line(-1,-2){3}} \put(82.3,5.0){\line(1,-2){3}}
	\put(26.4,5.4){\line(2,-3){3}} \put(103.6,5.4){\line(-2,-3){3}}
	\put(16.7,7.9){\line(-2,-3){3}} \put(113.3,7.9){\line(2,-3){3}}
	\put(8.7,27.5){\line(-1,-4){3}} \put(121.3,27.5){\line(1,-4){3}}
	\put(49.1,3.0){$a_2$} \put(55.0,8.9){$a_1$} \put(63.0,11.0){$a_0$} \put(71.0,8.9){$a_{n-1}$} \put(76.9,3.0){$a_{n-2}$}
	\put(34.4,11.5){$b_2$} \put(46.5,23.6){$b_1$} \put(63.0,28.0){$b_0$} \put(79.5,23.6){$b_{n-1}$} \put(91.6,11.5){$b_{n-2}$}
	\put(23.9,7.5){$c_2$} \put(35.5,25.8){$c_1$} \put(54.7,36.0){$c_0$} \put(76.5,35.1){$c_{n-1}$} \put(94.8,23.5){$c_{n-2}$} \put(105.0,4.3){$c_{n-3}$}
	\put(12.8,8.5){$d_2$} \put(26.2,31.8){$d_1$} \put(49.5,45.2){$d_0$} \put(76.5,45.2){$d_{n-1}$} \put(99.8,31.8){$d_{n-2}$} \put(113.2,8.5){$d_{n-3}$}
	\put(15.8,25.4){$e_2$} \put(37.5,45.4){$e_1$} \put(66.2,52.0){$e_0$} \put(94.4,43.2){$e_{n-1}$} \put(114.4,21.5){$e_{n-2}$}
	\put(6.0,28.5){$f_2$} \put(30.5,53.0){$f_1$} \put(64.0,62.0){$f_0$} \put(97.5,53.0){$f_{n-1}$} \put(122.0,28.5){$f_{n-2}$}
\end{picture}

	\caption{The graph of convex polytope $R_n''$.} \label{fig:rnsec}
\end{figure}

The graph of convex polytope $R_n'', n \geqslant 5$, is provided in Figure~\ref{fig:rnsec}.

\begin{Proposition}	For $n\geqslant 5$, it holds
	$$\gamma_{2-SRD}(R^{''}_n) =\left\lceil \frac{2|V(R^{''}_n)|}{3}\right\rceil=  4n.$$
\end{Proposition}
\proof

We define $g \colon V(R_n'') \mapsto \{0, 1,2\}$ as
\begin{align*}
	&g(a_i)  = g(c_i) = g(d_i)= g(f_i)=0, g(b_i) = g(e_i) = 2,\ i \in\{ 0, \ldots, n-1 \}.
\end{align*}

\begin{table}[H]
	\centering
	\begin{tabular}{c|l |c } \hline
		$v \in V_0$ & $N(v) \cap V_2$  & Indices  \\ \hline
		$a_i$      & $b_i$  &  \\
		$c_i$      &  $b_i$, $b_{i+1}$   & $i=0, \ldots, n-1$ \\
		$d_i$      &  $e_i$, $e_{i+1}$   &  \\
		$f_i$      &  $e_i$                &  \\ \hline
	\end{tabular}
	\caption{Coverage of  $V_0$ vertices of function $f$ defined for $R_n''$. } \label{opt-srdf-r_n_sec}
	
\end{table}

One can argue that $g$ is indeed a 2-SRD function following the arguments  given in Table~\ref{opt-srdf-r_n_sec}:
\begin{itemize}
	\item  Each vertex from $V_0$ is protected by at least one vertex from $V_2$;
	\item Each vertex from $V_2$ solely protects only one vertex from $V_0$. More precisely, each $b_i$ (respectively $e_i$) vertex protects solely  his neighboring $a_i$ (respectively $f_i$) vertex.
	
\end{itemize}

Further,
$$\omega(g) = 0 + 2n + 0 + 0 + 2n + 0  = 4n. $$

Therefore, $\gamma_{2-SRD}(R_n'') \leqslant 4n$. From this, and from Theorem~\ref{thm:cubic-lower-bound}, we obtain $\gamma_{2-SRD}(R_n'') = \left \lceil \frac{2}{3} |V(R_n'')|  \right \rceil = 4n. $ \\ \vspace{0.2cm}
\qed

	Convex polytope $A_n=(V(A_n),E(A_n)),n \geqslant 5$, or antiprism, was introduced in~\cite{baca1988labelings}. It is a graph with the following sets of vertices
	and edges
	\begin{align*}
		&V(A_n) = \{a_i, b_i, c_i \mid i \in \{ 0, \ldots, n - 1\} \},\\
		&E(A_n) = \{a_i a_{i+1}, b_ib_{i+1}, c_ic_{i+1}, a_ib_i, b_ic_i, a_{i+1}b_i, b_{i+1}c_i \mid i \in \{ 0,\ldots,  n-1\} \}.
	\end{align*}

	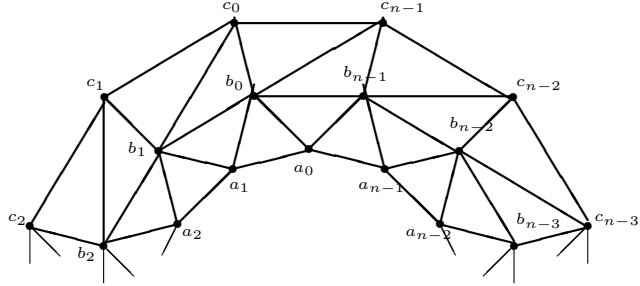
\begin{figure}
	\setlength\unitlength{1mm}
	\begin{picture}(130,60)
		\thicklines
		\tiny
		\put(47.7,5.0){\circle*{1}} \put(47.7,5.0){\line(1,1){7.3}} \put(47.7,5.0){\line(-4,-1){9.7}} \put(47.7,5.0){\line(-1,4){2.5}}
		\put(55.0,12.3){\circle*{1}} \put(55.0,12.3){\line(4,1){10.0}} \put(55.0,12.3){\line(-4,1){9.8}} \put(55.0,12.3){\line(1,4){2.8}}
		\put(65.0,15.0){\circle*{1}} \put(65.0,15.0){\line(4,-1){10.0}} \put(65.0,15.0){\line(-1,1){7.2}} \put(65.0,15.0){\line(1,1){7.2}}
		\put(75.0,12.3){\circle*{1}} \put(75.0,12.3){\line(1,-1){7.3}} \put(75.0,12.3){\line(-1,4){2.8}} \put(75.0,12.3){\line(4,1){9.8}}
		\put(82.3,5.0){\circle*{1}} \put(82.3,5.0){\line(1,4){2.5}} \put(82.3,5.0){\line(4,-1){9.7}}
		\put(38.0,2.2){\circle*{1}} \put(38.0,2.2){\line(3,5){7.2}} \put(38.0,2.2){\line(0,1){19.6}} \put(38.0,2.2){\line(-4,1){9.7}}
		\put(45.2,14.8){\circle*{1}} \put(45.2,14.8){\line(5,3){12.6}} \put(45.2,14.8){\line(3,5){10.0}} \put(45.2,14.8){\line(-1,1){7.1}}
		\put(57.8,22.0){\circle*{1}} \put(57.8,22.0){\line(1,0){14.5}} \put(57.8,22.0){\line(5,3){17.1}} \put(57.8,22.0){\line(-1,4){2.6}}
		\put(72.2,22.0){\circle*{1}} \put(72.2,22.0){\line(5,-3){12.6}} \put(72.2,22.0){\line(1,0){19.6}} \put(72.2,22.0){\line(1,4){2.6}}
		\put(84.8,14.8){\circle*{1}} \put(84.8,14.8){\line(3,-5){7.2}} \put(84.8,14.8){\line(5,-3){16.9}} \put(84.8,14.8){\line(1,1){7.1}}
		\put(92.0,2.2){\circle*{1}} \put(92.0,2.2){\line(4,1){9.7}}
		\put(28.3,4.8){\circle*{1}} \put(28.3,4.8){\line(3,5){9.8}}
		\put(38.1,21.9){\circle*{1}} \put(38.1,21.9){\line(5,3){17.0}}
		\put(55.2,31.7){\circle*{1}} \put(55.2,31.7){\line(1,0){19.7}}
		\put(74.8,31.7){\circle*{1}} \put(74.8,31.7){\line(5,-3){17.0}}
		\put(91.9,21.9){\circle*{1}} \put(91.9,21.9){\line(3,-5){9.8}}
		\put(101.7,4.8){\circle*{1}}
		\thinlines
		\put(47.7,5.0){\line(-1,-2){2}} \put(82.3,5.0){\line(1,-2){2}}
		\put(38.0,2.2){\line(1,-1){4}} \put(92.0,2.2){\line(-1,-1){4}}
		\put(38.0,2.2){\line(0,-1){5}} \put(92.0,2.2){\line(0,-1){5}}
		\put(28.3,4.8){\line(0,-1){5}} \put(101.7,4.8){\line(0,-1){5}}
		\put(28.3,4.8){\line(1,-1){4}} \put(101.7,4.8){\line(-1,-1){4}}
		\put(48.3,3.5){$a_2$} \put(54.5,9.7){$a_1$} \put(63.0,12.0){$a_0$} \put(71.5,9.7){$a_{n-1}$} \put(77.7,3.5){$a_{n-2}$}
		\put(34.5,0.5){$b_2$} \put(41.2,14.5){$b_1$} \put(54.0,23.3){$b_0$} \put(69.5,24.5){$b_{n-1}$} \put(83.5,17.8){$b_{n-2}$} \put(92.3,5.0){$b_{n-3}$}
		\put(25.4,5.4){$c_2$} \put(35.7,23.3){$c_1$} \put(53.6,33.6){$c_0$} \put(74.4,33.6){$c_{n-1}$} \put(92.3,23.3){$c_{n-2}$} \put(102.6,5.4){$c_{n-3}$}
		
		\end{picture}
		\caption{Convex polytope $A_n$}
	\end{figure}

	\Proposition For graph $A_n$, $n \geqslant 5$, it holds
	$$ \gamma_{2-SRD}(A_n)  \leqslant \begin{cases}
		\frac{3n}{2}, \textrm{if } n\equiv 0\bmod{2}, \\
		\lfloor \frac{3n}{2}\rfloor +1,  \textrm{otherwise}
		
	\end{cases}$$
	
	\proof

    We define the function $f=(V_0,V_1,V_2,V_3)$ according to the following cases:
	\begin{itemize}
	\item $ n\equiv 0\bmod{2}$;
	
	In this case, we set $V_3= \{b_{2i}: i\in\{1,\ldots, \frac{n}{2}\}\}$, $V_1=V_2=\emptyset$ and  $V_0=V\setminus V_3$.
	One can easily see that this yields a proper 2-SRD function $f$, as
	vertices $b_{2i}$ are capable of defending any two attacked vertices from $V_0$.
	We get $\omega(f)=3|V_3|=\frac{3n}{2}$.

	\item $ n\equiv 1\bmod{2}$;
	
	In this case, we set $V_3= \{b_{2i}: i\in\{1,\ldots, \frac{n-1}{2}\}\}$, $V_1=\{a_{n-1},c_{n-1}\}$, $V_2=\emptyset$ and
	$V_0=V\setminus (V_1\cup V_3)$.
	Analogously to the previous case,
	vertices $b_{2i}$ are capable of defending any two attacked  vertices  from $V_0$. Also,
	vertices $a_{n-1}$ and $c_{n-1}$ are labeled with 1, hence they are protected by themselves.
	So, function $f$ obtained in this way is a proper 2-SRD function with $\omega(f) = 3|V_3|+|V_1| = 3\frac{n-1}{2}+2 = \frac{3n+1}2=\lfloor\frac{3n}{2}\rfloor+1$.
    \end{itemize}
	\qed

 \subsection{Results for general $k$-SRD problem}\label{sec:general-graph-ksrd}
 At the beginning of this section, we provide a relation between the $k$-SRD number $\gamma_{k-\text{SRD}}$ and the Roman $k$-domination number  $\gamma_{kR}$.
 	
 	     \Proposition \label{proposition1} For arbitrary graph $G$, it holds $$\gamma_{k-\text{SRD}}(G)\leqslant  \gamma_{kR} (G).$$
 	
 	     \proof Let $f=(V_0,V_1,V_2)$ be a Roman $k$-dominating function, and $P$ an arbitrary attack pattern, containing $0<\sigma\leqslant k$ vertices from $V_0$. Let us denote these vertices with $v_1,v_2,\ldots ,v_{\sigma}$.  Following the definition of Roman $k$-dominating function, for each vertex $v_i$, $i\in\{1, \dotsc, \sigma\}$, there exists the set $S_{v_i}\subset V_2$ vertices which are adjacent with $v_i$ and $|S_{v_i}|\geqslant k$. We will explain the defence strategy for the pattern attack $P$. For  vertex $v_1$, we choose an arbitrary vertex $s_1\in S_{v_1}$, and for each other vertex $v_j$, $j\in\{2,\dotsc,\sigma\}$, we choose a vertex $s_j\in S_{v_j}\setminus \{s_1,s_2,\ldots,s_{j-1}\}$. Since the cardinality of each set $S_{v_i}$ is at least $k$, such a construction is clearly ensured.
 	   The chosen vertices $s_1,s_2, \ldots ,s_{\sigma}$ protect the vertices  $v_1,v_2,\ldots ,v_{\sigma}$ from the attack pattern $P$. Therefore, function $f$ satisfies conditions $(i)$ and $(ii)$ from Definition~\ref{def:SRDfunction}, which means $f$ is a proper $k$-SRD function. As $f$ is arbitrarily chosen, we conclude $$\gamma_{k-\text{SRD}}(G)\leqslant  \gamma_{kR} (G).$$
 	     \qed
 	
\noindent Another noticeable domination number known from literature is $k$-domination number $\gamma_{k}$. For graph $G=(V(G),E(G))$, the set $D \subset V(G)$ is a $k$-dominating set, if for every vertex $v \in V(G)$ it holds that $|N(v) \cap D| \geqslant k$. The $k$-domination number of graph $G$ is defined by  $\gamma_{k}(G) := \min  \{|D| \mid D \text{ is a proper } k\text{-dominating set on graph } G  \}$.
 In literature, a relation between  Roman $k$-domination number  $\gamma_{kR}$ and $k$-domination number $\gamma_{k}$ is established in~\cite{kammerling2009roman}, where inequality $ \gamma_{kR} (G)\leqslant 2 \gamma_{k} (G)$  for arbitrary graph $G$ is proven. Using this fact and Proposition~\ref{proposition1}, one may conclude

 	     \Corollary For arbitrary graph $G$, it holds 	
 	     $$\gamma_{k-\text{SRD}}(G)\leqslant  \gamma_{kR} (G)\leqslant 2 \gamma_{k} (G).$$   	
 	
\normalfont 	
Another upper bound can be derived from the basic domination number of graph $G$, denoted by $\gamma(G) = \gamma_1(G)$.

\Proposition  For any graph $G$ it holds $$  \gamma_{k-\text{SRD}} (G)\leqslant (k+1)\gamma(G).$$
\proof The claim follows from the fact that each vertex labelled with $k+1$ can protect any $k$ adjacent vertices. Therefore, if $S\subset V(G)$ is a domination set, it is sufficient to label each vertex from $S$ with $k+1$, and the remaining vertices with 0.
\qed \\

\Lemma For $k_1<k_2$, it holds $\gamma_{k_1-SRD}(G) \leqslant \gamma_{k_2-SRD}(G) \leqslant |V(G)|$.
\proof It follows directly from the fact that every feasible defense strategy against attack pattern $P$, with $|P|=k_2$, will ensure the protection against any attack pattern $P'\subset P$, with $|P'|=k_1$.  
\qed \\
\Lemma \label{lemma:2} For graphs $G = (V(G), E(G))$ and $H = (V(G), E(H))$ such that $E(G) \subset E(H)$, it holds  $\gamma_{k-\text{SRD}}(H) \leqslant \gamma_{k-\text{SRD}}(G)$.

\proof

Let $f$ be an arbitrary proper $k$-SRD function for the graph $G$. This function provides vertex labeling which determines a  feasible defense strategy
against $k$ attacks on graph $G$. Since graph $H$ has the same set of vertices as graph $G$ and preserves all existing connections between vertices of graph $G$, the specified defense strategy remains valid for graph $H$. In other words, function $f$ is also a proper $k$-SRD function for the graph $H$. As a consequence, we have $\gamma_{k-\text{SRD}}(H)\leqslant \omega(f)$. This inequality holds for an arbitrary $k$-SRD function $f$ for the graph $G$, so, in particular, we have $\gamma_{k-\text{SRD}}(H) \leqslant \gamma_{k-\text{SRD}}(G)$.
\qed \\

 Let $K_n =(V(K_n), E(K_n)), n \in \mathbb{N}$, be a complete graph, with $V(K_n)=\{v_i \mid i  \in \{1, \ldots, n\}\}$ and $E(V_n) =  \{ v_i v_j \mid i > j, i,j \in \{1, \ldots, n\}\}.$
 	\Lemma  \label{Kn_theorem} For complete graph $K_n$, it holds:
 	
 	  $$\gamma_{k-\text{SRD}}(K_n) =\begin{cases}
 		k+1, \textrm{if } k<n \\
 		n, \text{if } k=n.
 	\end{cases} $$
 	\proof
 	
 	 \underline{Upper bound}. Suppose $k < n$ and define a function $f \colon V(K_n) \mapsto \{0, \ldots, k+1\}$ by
 	 \begin{align*}
 	 	f(v_1) = k +1, f(v_i) = 0, i\in \{2, \ldots, n\}.
 	 \end{align*}
 	
 	 One can easily see that this mapping gives a proper $k$-SRD function as
 	 vertex $v_1$ is capable of defending any other  (attacked) $k$ vertices of graph $K_n$. Thus, $\gamma_{k-\text{SRD}}(K_n) \leqslant k+1$.
 	
 	  Suppose $k = n$.
 	  We can define a function $f \colon V(K_n) \mapsto \{0, \ldots, k\}$ by
 	 \begin{align*}
 	 	f(v_i) = 1, i\in \{1, \ldots, n\}.
 	 \end{align*}
 	
 	 Obviously, $f$ is a proper $k$-SRD function as
 	 every vertex is defended by itself. Therefore, $\gamma_{k-\text{SRD}}(K_n) \leqslant n$.
 	
 	 \underline{Lower bound}. Let $f$ be an arbitrary $\gamma_{k-\text{SRD}}$ function. For the sake of simplicity, we denote $|V_0| = \sigma$, and it is obvious that $\sigma <n$.
 	
 	 Suppose $k < n$.  We consider the following cases:
 	
 	 \begin{itemize}
 	 	\item $0 \leqslant \sigma < k$
 	 	
 	 	In this case, there are  $|V_{\geqslant 1}| = n-\sigma$ vertices labeled by at least 1, and they must protect $\sigma$ vertices labeled by zero, implying $\omega(f) \geqslant n-\sigma+\sigma = n \geqslant k+1$.
 	 	
 	 	\item $\sigma \geqslant k$
 	 	
 	 		Again, there are  $|V_{\geqslant 1}| = n-\sigma$ vertices labeled by at least 1, and some of them must protect $k$ vertices labeled by zero, implying $\omega(f) \geqslant n-\sigma+k \geqslant k+1$.
 	 \end{itemize}
 	
 	  Suppose $k = n$. Now, the only possible case is $0 \leqslant \sigma < k$, as $\sigma \geqslant k$ is not possible since $\sigma < n = k$.
 Similarly as in the previous case, there are  $|V_{\geqslant 1}| = n-\sigma$ vertices labeled with at least 1 which must protect $\sigma$ vertices labelled by zero, so $\omega(f) \geqslant n-\sigma+\sigma = n$. 
 	
     \qed \\
 	\Theorem \label{thm:extremal-upper-bound_n_eq_k}
 	Let $G=(V(G),E(G))$ be an arbitrary graph and $|V(G)|=n$.
 	
 	\begin{enumerate}
 		\item[(i)] 	If $k=n$, then $\gamma_{k-\text{SRD}}(G) = n$.
 		\item[(ii)] If $k<n$, then $\gamma_{k-\text{SRD}}(G) \geqslant k+1$.
 	\end{enumerate}

 	\proof
    Let $K_n=(V(G),E(K_n))$ be a complete graph. Since $E(G)\subseteq E(K_n)$, from Lemma \ref{lemma:2} follows $\gamma_{k-\text{SRD}}(K_n)\leqslant \gamma_{k-\text{SRD}}(G)$,
    so, by applying Lemma \ref{Kn_theorem}, we obtain  $\gamma_{k-\text{SRD}}(G)\geqslant \begin{cases}
 		k+1, \textrm{if } k<n \\
 		n, \text{if } k=n
 	\end{cases}.$ For $k=n$, all vertices from $V(G)$ can be labeled with $1$, implying that, in this case, inequality $\gamma_{k-\text{SRD}}(G)\leqslant n$ also holds.
 \qed \\
%
%
%

\Corollary  \label{cor:extremal-vertex_degree_n-1}	Let $G=(V(G),E(G))$ be an arbitrary graph  and $|V(G)|=n$. If there exists a vertex $v\in V(G)$ such that $d(v)=n-1$, then
 	  $$\gamma_{k-\text{SRD}}(G) =\begin{cases}
	k+1, \textrm{if } k<n \\
	n, \text{if } k=n.
\end{cases} $$

\normalfont
	Let $P_n=(V(P_n),E(P_n))$ be a path graph, with $V(P_n)=\{a_i:i\in\{1,2,\dotsc,n\}\}$ and $E(P_n)=\{a_ia_{i+1}:i\in\{1,2,\dotsc,n-1\}\}$.

    \Theorem For  $k \geqslant 2$, it holds $\gamma_{k-\text{SRD}}(P_n) =
    \left\lceil  \frac{2k}{2k+1} n \right\rceil$.
    \proof

  Since $\Delta({P_n})=2$, every $k$-SRD function on $P_n$ can be written as $f=(V_0,V_1,V_2,V_3)$.  At this moment, we state an important (and obvious) observation concerning the relation between sets $V_0$ and $V_2$: Any $k$ different vertices from $V_0$ cannot be covered solely by $k-1$ or fewer vertices from $V_2$; however, any $k+1$ different vertices from $V_0$ can be covered solely by $k$ vertices from $V_2$.

    \underline{Upper bound}.     Let us define the function $f:V(P_n)\rightarrow\{0,1,2,3\}$ by:

    \begin{itemize}
    	\item For $n\in\{k,k+1,\dotsc,2k\}$, we set $f(a_i):=1,i\in\{1,2,\dotsc, n\}$.
    	\item For $n=(2k+1)l+r$, with $l\geqslant 1$, $r\in\{0,1,\dotsc,2k\}$, we introduce the block $B:=02020\underbrace{\underline{20}\;\underline{20}\dotsc \underline{20}}_{k-2}$
    	of the length $2k+1$, as the pattern for defining values for $2k+1$ consecutive vertices from $V(P_n)$. More precisely, for an arbitrary $j\in\{1,2,\dotsc, l\}$,
    	we set the values
    	$$f(a_{(2k+1)(j-1)+1}),f(a_{(2k+1)(j-1)+2}),\dotsc,f(a_{(2k+1)j}),$$
    	to be equal to the values located at the same positions in the block $B$. In other words, the values of function $f$ for the first $(2k+1)l$ vertices are obtained
    	by concatenating $l$ copies of the block $B$. Besides that, values of all remaining vertices from the ``tail'' of
    	the path are labelled with one, i.e., $f(a_{(2k+1)l+s}):=1$, $s\in\{1,2\dotsc,r\}$.
    \end{itemize}

    Inside the block $B$ there are a total of $k+1$ zeros and they are covered by $k$ twos. Also, any different $m\leqslant k$ zeros are covered by exactly $m$ twos. In accordance with our
  previously stated observation, this results in a valid defence configuration yielding a proper $k$-SRD function. Furthermore, for $l\geqslant 1$ and $r\in\{0,1,\dotsc,2k\}$, it holds

    $$\omega(f)=\left\{
    \begin{array}{ll}
    	n, & \hbox{for $n\in\{k,k+1,\dotsc,2k\}$;} \\
    	2kl+r, & \hbox{for $n=(2k+1)l+r$;}
    \end{array}
    \right.$$
    For $n\in\{k,k+1,\dotsc,2k\}$, we have $n = \left\lceil n-\frac{n}{2k+1}\right\rceil = \left\lceil\frac{2k}{2k+1}n\right\rceil$. Also,
    for $n=(2k+1)l+r$, with $l\geqslant 1$ and $r\in\{0,1,\dotsc,2k\}$, it holds
    $$ 2kl+r = 2kl+\left\lceil r-\frac{r}{2k+1}\right\rceil = \left\lceil 2kl+\frac{2kr}{2k+1}\right\rceil =\left\lceil\frac{2k}{2k+1}n\right\rceil.  $$
    Putting all these facts together, we conclude $ \gamma_{k-\text{SRD}}(P_{n})\leqslant \left\lceil \frac{2k}{2k+1}n\right\rceil. $

    \underline{Lower bound}. Let us assume that $g=(V_0,V_1,V_2,V_3)$ is an arbitrary $k$-SRD function defined on $V(P_n)$. 
  We need to prove inequality $\omega(g)\geqslant \left\lceil \frac{2k}{2k+1}n\right\rceil$. To do so, it is sufficient to describe the appropriate discharging procedure in which weights of vertices from $V_2\cup V_3$ are used in order to charge weights of vertices from $V_0$ and accomplishing that, in the end, every vertex from  $V(P_n)$ would have a weight of at least $\frac{2k}{2k+1}$.

  First, note that every vertex $a_i\in V_0$ such that $N(a_i)\cap V_3\ne\emptyset$ can receive the full charge of $\frac{2k}{2k+1}$
  from a neighbor in $N(a_i)\cap V_3$. Moreover, every vertex $a_j\in V_3$ such that $N(a_j)\subseteq V_0$ can charge all of his neighbors with full charge of
  $\frac{2k}{2k+1}$ each, since after this discharge, the remaining weight of this vertex is $3-\frac{4k}{2k+1}=\frac{2k+3}{2k+1}>1>\frac{2k}{2k+1}$. So, w.l.o.g. we assume $g=(V_0,V_1,V_2)$ and describe the discharging procedure related to this type of $k$-SRD function.


   Let us observe an arbitrary subpath of length $2m+1$, $m<k$, which contains $m$ zeros. Since those zeros cannot be covered by $m-1$ or less twos, we conclude that,
   within the observed path, every zero has a unique two as its defender. This ``1--1''  correspondence allows us to charge every zero within this subpath with $\frac{2k}{2k+1}$ from its unique defender. The remaining weight of each defender within this block is $2-\frac{2k}{2k+1}=\frac{2k+2}{2k+1}>1>\frac{2k}{2k+1}$.

\begin{figure}
	
    \adjustbox{scale=0.75,center}{
    	\begin{tikzcd}
    		0  & 2 \arrow[l,dashrightarrow, bend left=50,"{\frac{2k}{2k+1}}"] \arrow[r,dashrightarrow, swap, bend right=50,"{\frac{2}{2k+1}}"]& 0 &
    		2  \arrow[l,dashrightarrow, bend left=50,"{\frac{2k-2}{2k+1}}"] \arrow[r,dashrightarrow, swap, bend right=50,"{\frac{4}{2k+1}}"] & 0 & \dotsc & 0 &
    		2  \arrow[l,dashrightarrow, bend left=50,"{\frac{4}{2k+1}}"] \arrow[r,dashrightarrow, swap, bend right=50,"{\frac{2k-2}{2k+1}}"]   & 0 &
    		2  \arrow[l,dashrightarrow, bend left=50,"{\frac{2}{2k+1}}"] \arrow[r,dashrightarrow, swap, bend right=50,"{\frac{2k}{2k+1}}"] &  0
    	\end{tikzcd}

    }

\caption{The memory-based discharging  for $P_n$: the case of an  existence of $2k+1$ subpath with $k+1$ zeros.}\label{fig:memory-dischargind-case-3}
\end{figure}
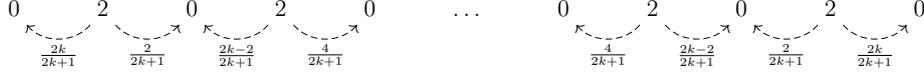

   Finally, the most challenging case is the existence of a subpath of length $2k+1$ which contains $k+1$ zeros. Then, it is possible to cover those zeros with $k$ twos, by alternating zeros and twos (basically, this is exactly block $B$  introduced in the proof of the upper bound). The memory-based discharging procedure which are we looking for is shown in
   Figure~\ref{fig:memory-dischargind-case-3}. 

    More in details, the most left vertex labelled with a two gives  $\frac{2k}{2k+1}$ charge to its left zero neighbor, and $\frac{2}{2k+1}$ charge to its right zero neighbor, the second from left vertex labelled with a
   two  gives  $\frac{2k-2}{2k+1}$ charge to its left zero neighbor, and  $\frac{4}{2k+1}$ charge to its right zero neighbor, and, generally, the $m-$th
   leftmost two, $m\in\{1,2,\dotsc, k\}$, give $\frac{2k-2(m-1)}{2k+1}$ charge to its left zero neighbor, and $\frac{2m}{2k+1}$ charge to its right zero neighbor.
   Clearly, after this charging procedure, every zero is charged by $\frac{2k}{2k+1}$. Also, for an arbitrary $m\in\{1,2,\dotsc,k\}$, the remaining weight of every
   defender is $2-\frac{2k+2}{2k+1}=\frac{2k}{2k+1}$. The depicted discharging procedure can be generalized to the case of a subpath of length $2s+1$,$s>k$, which
   contains $s+1$ zeros covered by $s$ twos. Simply, from this subpath we can extract a block of length $2k+1$ which contains $k+1$ zeros covered with $k$ twos. After discharging
   this block in the described manner, this procedure can be repeated until all ``leftovers''  after some iteration do not become subpaths with lengths smaller than $2k+1$. Then, we exploit the existence of ``1--1'' correspondence between zeros and twos to finish off the discharging procedure.

    In each of these cases, we ensured the weights of all vertices of graph $P_n$ are greater or equal to $\frac{2k}{2k+1}$, which implies $\omega(g)\geqslant \frac{2k}{2k+1}n$. Thus, we have $\gamma_{k-\text{SRD}}(P_n)\geqslant \frac{2k}{2k+1}n$, and since $\gamma_{k-\text{SRD}}$ is an integer value, we conclude
    $ \gamma_{k-\text{SRD}}(P_{n})\geqslant \left\lceil \frac{2k}{2k+1}n\right\rceil. $ \\
    \qed \\
    Let $C_n=(V(C_n),E(C_n))$ be a cycle graph, with $V(C_n) = \{a_i \colon i \in \{1, \ldots, n \}\}$ and $E(C_n)= \{ a_i a_{i+1} \colon i \in \{1, \ldots, n-1\}\} \cup \{a_na_1\}$.

  \Proposition  For  $k \geqslant 2$, it holds $\gamma_{k-\text{SRD}}(C_n) =
  \left\lceil  \frac{2k}{2k+1} n \right\rceil$.

  \proof
    The proof relies on the previously established facts. Indeed, we can perform ``unchaining'' of graph $C_n$ to get a path graph $P_n=(V(P_n),E(P_n))$, such that $V(P_n):=V(C_n)$ and $E(P_n):=E(C_n)\setminus\{a_na_1\}$. Now, using
    Lemma \ref{lemma:2}, we get $\gamma_{k-\text{SRD}}(C_n)\leqslant \gamma_{k-\text{SRD}}(P_n)=\left\lceil \frac{2k}{2k+1}n\right\rceil.$ On the other hand, for an arbitrary $\gamma_{k-\text{SRD}}$ function $f$ on graph $C_n$, we have the following two cases:
    \begin{itemize}
    	\item $f(a_i)>0$, for all $i\in\{1,2,\ldots,n\}$;
    	
    	In this case, by removing any edge from $E(C_n)$, function $f$ would still be a  $\gamma_{k-\text{SRD}}$ function on the  obtained path graph.
    	\item $f(a_{i_0})=0$, for some $i_0\in\{1,2,\dotsc, n\}$;
    	
    	In this case, vertex $a_{i_0}$ has a neighbor labeled with at least two. Deleting the edge not leading to that neighbor, we obtain a path graph
    	on which function $f$ is again a proper $k$-SRD function.
    \end{itemize}
    So, in both cases we get $\omega(f)\geqslant \frac{2k}{2k+1}n$, and consequently $\gamma_{k-\text{SRD}}(C_n)\geqslant \frac{2k}{2k+1}n$. Since
    $\gamma_{k-\text{SRD}}$ is an integer value, it follows $\gamma_{k-\text{SRD}}(C_n)\geqslant\left\lceil \frac{2k}{2k+1}n\right\rceil.$

    Therefore, $\gamma_{k-\text{SRD}}(C_n)=\left\lceil \frac{2k}{2k+1}n\right\rceil$. 
    \qed \\ \vspace{0.2cm}

    Let $S_n=(V(S_n),E(S_n))$ be a star graph with $V(S_n)= \{a_i \colon i \in \{1, \ldots, n\}\}$ and $E(S_n)= \{ a_1 a_i \colon  i \in \{ 2, \ldots, n  \}$. Wheel graph $W_n=(V(W_n), E(W_n))$ is a graph with $V(W_n)= \{a_i \colon i \in \{1, \ldots, n\}\}$ and $E(W_n) = \{ a_1 a_i \colon i \in \{2, \ldots, n \}\}  \cup \{ a_i a_{i+1} \colon i \in \{1, \ldots, n-1 \}\} \cup \{a_n a_1\}$.

    \Proposition \label{thm:s_n} For $k\geqslant 2$ and $n \geqslant 3$,  it holds

    $$\gamma_{k-\text{SRD}}(S_n)
        =\begin{cases}
    	    k+1, \textrm{if } n>k \\
    	    n, \text{if } n=k.
    \end{cases} $$

     The same exact value applies for wheel graphs $W_n, n \geqslant 4$.
  \proof     The claim is a direct consequence of Corollary~\ref{cor:extremal-vertex_degree_n-1}.
\vspace{0.2cm}

Let $V_u = \{u_1, \ldots, u_n\}$ and $V_v = \{ v_1, \ldots, v_p\}$, $n \geqslant p \geqslant 1$.
Graph $K_{n, p} = (V(K_{n,p}), E(K_{n,p}))$, $n \geqslant p \geqslant 1$, is called complete bipartite graph if $V(K_{n,p})=V_u\cup V_v$, $V_u \cap V_v = \emptyset$,
and $E(K_{n,p}) = \{ u_i v_j \mid i \in \{1, \ldots, n\}, j \in \{1, \ldots, p\} \}$. Note that $K_{n,1}$ graph is essentially a star graph $S_{n+1}$. \\


Prior to formulating the theorem pertaining to complete bipartite graphs, we introduce requisite notation and present important observations used in the proof. Consider an arbitrary function  $k$-SRD function $f$  operating on a complete bipartite graph $K_{n, p}$, where $n \geqslant p \geqslant 1$. The number of vertices labeled with 0 in the partitions $V_u$ and $V_v$ are denoted as $\sigma_u$ and $\sigma_v$ respectively. Let $a$ denote either of the indices $u$ or $v$, and let $b$ represent the other one.

In the subsequent observation, we undertake an analysis of the lower bound of $\gamma_{k-\text{SRD}}(K_{n, p})$, depending on  the count of zero-labeled vertices in the partitions $V_u$ and $V_v$.

\Observation \label{observationKNP}
The lower bounds of an arbitrary  $k$-SRD function $f$ on complete bipartite graph $K_{n, p}$ w.r.t.\  $\sigma_u$ and  $\sigma_v$, are as presented in the following cases:
\normalfont

\begin{enumerate}
	
\item [(i)] $\sigma_a = \sigma_b = 0$

This refers to the case when no vertices are labeled with zero. Here all $n+p$ vertices are  labeled with at least 1, thus  $\omega(f) \geqslant n+p$.

\item [(ii)] $\sigma_a =0$ and $0< \sigma_b <k$

In this case, we have:
\begin{itemize}
	\item all $|V_a|$  $a$-vertices are labeled with at least 1,
	\item some of $a$-vertices must protect $\sigma_b$ $b$-vertices labeled with 0,
	\item there are at least $|V_b|-\sigma_b$ $b$-vertices labeled with at least 1.
\end{itemize}

Therefore, $\omega(f) \geqslant |V_a|+\sigma_b+|V_b|-\sigma_b = |V_a|+|V_b|=n+p$.

\item [(iii)] $\sigma_a =0$ and $\sigma_b \geqslant k$

In this case, we have:
\begin{itemize}
	\item all $|V_a|$ $a$-vertices are labeled with at least 1,
	\item some of $a$-vertices must protect any $k$ $b$-vertices labeled with 0,
	\item there are at least $|V_b|-\sigma_b$ $b$-vertices labeled with at least 1.
\end{itemize}

Therefore, $\omega(f) \geqslant |V_a|+k+|V_b|-\sigma_b \geqslant |V_a|+k  \geqslant \min\{|V_a|,|V_b|\}+k = p+k$.

\item [(iv)] $0< \sigma_a <k$ and $ 0<\sigma_b <k$

In this case, we have:

\begin{itemize}
	\item there are at least $|V_a|-\sigma_a$ $a$-vertices labeled with at least 1,
	\item some of $a$-vertices must protect $\sigma_b$ $b$-vertices labeled with 0,
	\item there are at least $|V_b|-\sigma_b$ $b$-vertices labeled with at least 1,
	\item some of $b$-vertices must protect $\sigma_a$ $a$-vertices labeled with 0.
\end{itemize}

We conclude $\omega(f) \geqslant |V_a|-\sigma_a+\sigma_b+|V_b|-\sigma_b+\sigma_a = |V_a|+|V_b| =n+p$.

\item [(v)] $0< \sigma_a <k$ and $ \sigma_b \geqslant k$

In this case, we have:

\begin{itemize}
	\item there are at least $|V_a|-\sigma_a$ $a$-vertices labeled with at least 1,
	\item some of $a$-vertices must protect any $k$ $b$-vertices labeled with 0,
	\item there are at least $|V_b|-\sigma_b$ $b$-vertices labeled with at least 1,
	\item some of $b$-vertices must protect $\sigma_a$ $a$-vertices labeled with 0.
\end{itemize}

Thus, we obtain $\omega(f) \geqslant |V_a|-\sigma_a+k+|V_b|-\sigma_b+\sigma_a \geqslant |V_a| +k+1$, since  $\sigma_b<|V_b|$ because $\sigma_a>0$.

\item [(vi)] $\sigma_a \geqslant k$ and $\sigma_b \geqslant k$

In this case we have:
\begin{itemize}
	\item there are at least $|V_a|-\sigma_a$ $a$-vertices labeled with at least 1,
	\item some of $a$-vertices must protect any $k$ $b$-vertices labeled with 0,
	\item there are at least $|V_b|-\sigma_b$ $b$-vertices labeled with at least 1,
	\item some of $b$-vertices must protect any $k$ $a$-vertices labeled with 0.
\end{itemize}

 We obtain
$\omega(f) \geqslant |V_a|-\sigma_a + k + |V_b|-\sigma_b + k \geqslant 1+ k+1+k = 2k+2$, since $|V_a|>\sigma_a$ and $|V_b|>\sigma_b$.

\end{enumerate}
\qed

\Theorem For graph $K_{n, p}, n \geqslant p >1$, it holds

$$\gamma_{k-\text{SRD}}(K_{n, p}) =\begin{cases}
	2k+2, \textrm{if } k<p-1; \\
	\min\{n,k\}+p, \text{if } k\geqslant p-1.\\
\end{cases} $$

\proof

\underline{\textit{Upper bound}}.
Suppose $k < p-1$. We define a function $f:V(K_{n,p})\rightarrow\{0,1,\dotsc,k+1\}$ by
$$ f(u_1) = k+1,  f(v_1) = k+1, f(u_i) = 0, i \in \{2, \ldots, n\}, f(v_j) = 0, j \in \{2, \ldots, p\}.$$

One can easily verify that $f$ gives indeed a $k$-SRD function, $\omega(f) = 2k+2$, thus $\gamma_{k-\text{SRD}}(K_{n,p})\leq 2k+2$.

Suppose now that $k \geqslant p-1$. We define a function $f:V(K_{n,p})\rightarrow\{0,1,\dotsc,\min\{n,k\}+1\}$ by
$$ f(u_i) = 0, i \in \{1, \ldots, n\},  f(v_1) = \min\{n,k\}+1,  f(v_j) = 1, j \in \{2, \ldots, p\}.$$

Since vertex $v_1$ is capable to solely defend all $u$-vertices against any of $n\choose \min\{n,k\}$ attack patterns, it follows that $f$ is a proper $k$-SRD function
such that $\omega(f) = \min\{n,k\}+1 + (p-1) = \min\{n,k\}+p$, which implies $\gamma_{k-\text{SRD}}(K_{n,p}) \leq \min\{n,k\} + p$.

\underline{Lower bound}. Let $g$ be an arbitrary  $k$-SRD function defined on the vertices of graph $K_{n,p}$.  By utilizing Observation~\ref{observationKNP}, we cover all cases based on the number of $u$ and $v$ vertices labelled with zero. In Tables~\ref{tbl:lowerbounds1}--\ref{tbl:lowerbounds2} we show the lower bounds for cases $k < p-1$ and $k\geqslant p-1$, respectively. Both of these tables follow a similar organization. The first column contains possible ranges for the value $\sigma_u$, representing the number of $u$ vertices labelled with 0. Similarly, the second column shows this information for $v$ vertices. The third column indicates the case in Observation~\ref{observationKNP} that the combination of $\sigma_u$ and $\sigma_v$ falls into. The last column displays the calculation of the considered lower bound.

For instance, consider the second row of Table~\ref{tbl:lowerbounds1}, representing the case where there are no $u$ vertices labeled with zero ($\sigma_u = 0$ in the first
column), while there are some, but fewer than $k$, $v$ vertices labelled with zero  ($0<\sigma_v<k$ in the second column).
The third column indicates that this situation falls into case $(ii)$ of Observation~\ref{observationKNP}, where $a=u$ and $b=v$. The calculation  follows from  case $(ii)$ of Observation~\ref{observationKNP}, and relies on $n\geqslant p$ and  $k < p-1$.

From the last column of the tables, one can clearly notice that $\omega(g)$ is greater than or equal to the proposed lower bound for all cases, which concludes the proof.

\begin{table}[h]\small
	\caption{Lower bounds of $\omega(g)$, for the case $k < p-1$.}
	\begin{tabular}{|c|c|l|l|}\hline
		$\sigma_u$             & $\sigma_v$             & Observation case & Calculation                                                               \\\hline
		                    & 0                      & ($i$):   $a=u, b=v$          & $\omega(g) \geqslant n+p \geqslant 2p > 2k+2$                                          \\
		0                      & $0<\sigma_v<k$         & ($ii$): $a=u, b=v$             & $\omega(g) \geqslant n+p \geqslant 2p > 2k+2$                             \\
		                    & $\sigma_v \geqslant k$ & ($iii$): $a=u, b=v$            & $\omega(g) \geqslant p+k > 2k+1$, $\omega(g)  \geqslant 2k+2$  \\\hline
		      & 0                      & ($ii$): $b=u, a=v$         & $\omega(g) \geqslant n+p \geqslant 2p > 2k+2$                             \\
		$0<\sigma_u<k$         & $0<\sigma_v<k$         & ($iv$): $a=u, b=v$             & $\omega(g) \geqslant n+p \geqslant 2p > 2k+2$                             \\
	         & $\sigma_v \geqslant k$ & ($v$): $a=u, b=v$              & $\omega(g) \geqslant n+k+1 > 2k+2$                                \\\hline
		 & 0                      & ($iii$): $b=u, a=v$        & $\omega(g)  \geqslant   p+k > 2k+1$, $\omega(g)  \geqslant 2k+2$ \\
		$\sigma_u \geqslant k$ & $0<\sigma_v<k$         & ($v$): $b=u, a=v$         & $\omega(g)    \geqslant  p + k+1 >   2k+2$   \\
	 & $\sigma_v \geqslant k$ & ($vi$): $a=u, b=v$         & $\omega(g)  \geqslant   2k+2$            \\\hline
	\end{tabular}\label{tbl:lowerbounds1}
\end{table}

\begin{table}[h]\small
		\caption{Lower bounds of $\omega(g)$, for the case $k \geqslant p-1$.}
	\begin{tabular}{|c|c|l|p{6.5cm}|}\hline
		$\sigma_u$                                & $\sigma_v$             & Observation case    & Calculation\\\hline
		\multirow{3}{*}{0}                        & 0                      & ($i$):  $a=u, b=v$  & $\omega(g) \geqslant n+p \geqslant \min\{n,k\}+p$\\
		& $0<\sigma_v<k$         & ($ii$): $a=u, b=v$  & $\omega(g) \geqslant n+p    \geqslant \min\{n,k\}+p$ \\
		& $\sigma_v \geqslant k$ & ($iii$): $a=u, b=v$ & $\omega(g)  \geqslant   p+k \geqslant \min\{n,k\}+p$                          \\\hline
		\multirow{3}{*}{$0<\sigma_u<k$}           & 0                      & ($ii$): $b=u, a=v$  & $\omega(g) \geqslant n+p \geqslant \min\{n,k\}+p$                                \\
		& $0<\sigma_v<k$         & ($iv$): $a=u, b=v$  & $\omega(g) \geqslant    n+p \geqslant \min\{n,k\}+p$                             \\
		& $\sigma_v \geqslant k$ & ($v$): $a=u, b=v$   & $\omega(g)\geqslant    n+k+1 \geqslant n+p \geqslant \min\{n,k\}+p$              \\\hline
		\multirow{3}{*}{$\sigma_u \geqslant   k$} & 0                      & ($iii$): $b=u, a=v$ & $\omega(g)  \geqslant k+p   \geqslant \min\{n,k\}+p$                                      \\
		& $0<\sigma_v<k$         & ($v$): $b=u, a=v$   & $\omega(g) \geqslant p+k+1 > \min\{n,k\}+p$                                      \\
		& $\sigma_v \geqslant k$ & ($vi$): $a=u, b=v$  & $\omega(g)  \geqslant   2k+2 \geqslant \min\{n,k\} +1 + p-1 +1 > \min\{n,k\} +p$\\\hline
	\end{tabular}\label{tbl:lowerbounds2}
\end{table}
\qed

Let $P_2\times P_n=(V(P_2\times P_n),E(P_2\times P_n))$ be a grid graph, with $V(P_2\times P_n)=\{a_{ij}:i\in\{1,2\},j\in\{1,\dotsc,n\}\}$ and
$E(P_2\times P_n)=\{a_{ij}a_{ij+1}:i\in\{1,2\},j\in\{1,\dotsc,n-1\}\}\cup\{a_{1j}a_{2j}:j\in\{1,\dotsc,n\}\}$.
 	
\Theorem For an arbitrary $k\geqslant 2$ and grid graph $P_2 \times P_n$, $n \geqslant \frac{k}{2}$, it holds
$$ \gamma_{k-\text{SRD}}(P_2 \times P_n) = \begin{cases}
	   \lceil \frac{2kn}{k+1} \rceil, \textrm{if } n \equiv r \bmod (k+1) \textrm{ for some }  0<r<\frac{k+1}{2}; \\
	   \lceil \frac{2kn}{k+1} \rceil + 1, \mbox{otherwise}.
\end{cases}  $$

\proof Note that every proper $k$-SRD function of the grid graph $P_2\times P_n$ takes labels from $\{0,1,2,3\}$ in the case $k=2$ and labels from
 $\{0,1,2,3,4\}$ in the case $k>2$. In the rest of this proof, we denote the appropriate set of labels as $S$. Additionally, for a matrix $M$ of format $2\times n$, the matrices derived from $M$ by using upside-down reflection of its rows and left-to-right reflection of its columns will be denoted by $M^{\updownarrow}$ and $M^{\leftrightarrow}$, respectively.

\underline{Upper bound}. Let us define a function $f:V(P_2\times P_n)\rightarrow S$ in the following manner:
\begin{itemize}
  \item For $n\in\left\{\frac{k}{2},\dotsc,k\right\}$, we set $f(a_{ij}):=1,i\in\{1,2\}, j\in\{1,2,\dotsc, n\}$.
  \item For $n\in\left\{k+1,\dotsc,2k+1\right\}$, we construct a matrix $M_n$ of the format $2\times n$ with the idea to initialize values $f(a_{ij})$ to be equal to the corresponding values of this matrix. This construction is started by introducing the ``central'' matrix $$M_0:=\left (
\begin{array}{rrrrrrr}
0 & 2 & 0 & 2 & 0 & \dotsc &  \\
2 & 0 & \undermat{k-2}{2 & 0 & 2 & \dotsc}\\
\end{array}
\right )$$

of the format $2\times k$. Then this matrix is expanded with appropriate $n-k$ columns, such that $\left\lceil\frac{n-k}{2}\right\rceil$ columns are added to the right end side of the matrix $M_0$ and the rest of $\left\lfloor\frac{n-k}{2}\right\rfloor$ columns are added to the left end side of the matrix $M_0$. More specifically, for an integer $m\geqslant 1$, we define a matrix $B_m:=\left(
\begin{array}{rrrrrrr}
1 & 0 & 2 & 0 & 2 & 0 & \dotsc  \\
\undermat{m}{0 & 2 & 0 & 2 & 0 & 2 & \dotsc} \\
\end{array}
\right )$ of format $2\times m$, and set

$$M_n:=\left\{
         \begin{array}{ll}
           \left (
\begin{array}{r|r|r}
B_{\left\lfloor\frac{n-k}{2}\right\rfloor}^{\leftrightarrow} & M_0 & B_{\left\lceil\frac{n-k}{2}\right\rceil}^{\updownarrow} \\
\end{array}
\right ), & \hbox{for even $k$;} \\
\left (
\begin{array}{r|r|r}
B_{\left\lfloor\frac{n-k}{2}\right\rfloor}^{\leftrightarrow} & M_0 & B_{\left\lceil\frac{n-k}{2}\right\rceil} \\
\end{array}
\right ) , & \hbox{for odd $k$.}
         \end{array}
       \right.
$$

\item For $n=(k+1)l+r$, with $l\geqslant 2$ and $r\in\{0,1,\dotsc,k\}$, the goal is to construct a matrix $M_n$ of format $2\times n$ and initialize values $f(a_{ij})$ to be equal to the corresponding elements of this matrix. Since $n-(k+1)(l-1)=k+1+r\in\{k+1,\dotsc,2k+1\}$, we can use the construction from the previous case to obtain a matrix $M_{n-(k+1)(l-1)}$.  Consider a matrix $P_0$ of the format $2\times(k+1)$ defined by
$$P_0:=\left (
\begin{array}{rrrrr|rrr|rrrrr}
\dotsc & 2 & 0 & 2 & 0 & 2 & 0 & 0  & 2 & 0 & 2 & 0 &  \dotsc  \\
\undermat{\left\lceil\frac{k-2}{2}\right\rceil}{\dotsc & 0 & 2 & 0 & 2 } & 0 & 0 & 2 & \undermat{\left\lfloor\frac{k-2}{2}\right\rfloor}{ 0 & 2 & 0 & 2 &  \dotsc}\\
\end{array}
\right )$$

\vspace{10pt}
The next step includes creating a matrix $P_1$ of format $2\times (k+1)(l-1)$ by adjoining $l-1$ copies of the matrix $P_0$ in a way towards preserving feasibility of function $f$. This can be done by putting
$$P_1:=\left\{
         \begin{array}{lll}
           \left (
\begin{array}{r|r|r|r|r}
\undermat{(k+1)(l-1)}{P_0 & P_0 & P_0 & P_0 & \dotsc} \\
\end{array}
\right ), & \hbox{for even $k$;} \\
\\
\left (
\begin{array}{r|r|r|r|r}
\undermat{(k+1)(l-1)}{P_0^{\updownarrow} & P_0 & P_0^{\updownarrow} & P_0 & \dotsc} \\
\end{array}
\right ) , & \hbox{for odd $k$.}
         \end{array}
       \right.
$$

\vspace{5pt}

Finally, we insert pattern matrix $P_1$ between two consecutive columns of matrix $M_{n-(k+1)(l-1)}$. More precisely, we place $P_1$ into the middle region of the central submatrix $M_0$ of matrix $M_{n-(k+1)(l-1)}$ by inserting pattern matrix $P_1$ between two central columns of $M_0$ in the case of even $k$ or between middle column of $M_0$ and its predecessor column in the case of odd $k$. After that, optionally we apply an upside-down reflection of $P_1$, to ensure that the left jointed section consists of two different columns. At the end, if the right jointed section consists of two equal columns, we perform another upside-down reflection of the remaining right side submatrix of the matrix $M_{n-(k+1)(l-1)}$.

\end{itemize}

The previous construction is performed such that any $m\leqslant k$ different zeros in matrix $M_n,n\geqslant\frac{k}{2}$, are covered by at least $m$ different twos. Hence, the constructed function $f$ that corresponds to this matrix gives a proper $k$-SRD function of grid graph $P_2\times P_n$.

Clearly, for $n\in\left\{\frac{k}{2},\dotsc,k\right\}$, it holds $\omega(f)=2n$. The weights inside matrices $M_0$ and $P_0$ are both equal to $2k$ and the weight inside matrix $B_m$ is equal to $2m-1$. Consequently, for $n=(k+1)l,l\geqslant 1$, it holds
$\omega(f)=(2k+1)+2k(l-1)=2kl+1$, and for $n=(k+1)l+r,l\geqslant 1,r\in\{1,2,\dotsc,k\}$, it holds $\omega(f)=(2k+2r)+2k(l-1)=2kl+2r$. On the other hand, for
$n\in\left\{\frac{k}{2},\dotsc,k\right\}$, it holds
$$ \left\lceil \frac{2k}{k+1}n \right\rceil=\left\lceil 2n-\frac{2n}{k+1}\right\rceil=\left\{
                                                                                         \begin{array}{ll}
                                                                                          2n , & \hbox{for $n\leqslant\frac{k+1}{2}$;} \\
                                                                                          2n-1 , & \hbox{otherwise.}
                                                                                         \end{array}
                                                                                       \right.
$$
Also, for $n=(k+1)l,l\geqslant 1$, it holds $\left\lceil \frac{2k}{k+1}n \right\rceil=2kl$, while for $n=(k+1)l+r,l\geqslant 1,r\in\{1,2,\dotsc,k\}$, it holds
\begin{align*} \left\lceil \frac{2k}{k+1}n \right\rceil&=\left\lceil 2kl+\frac{2kr}{k+1}\right\rceil=2kl+\left\lceil 2r-\frac{2r}{k+1}\right\rceil\\
&=\left\{
                                                                                                                                        \begin{array}{ll}
                                                                                                                                         2kl+2r, & \hbox{for $0<r<\frac{k+1}{2}$;} \\
                                                                                                                                          2kl+2r-1, & \hbox{otherwise.}
                                                                                                                                        \end{array}
                                                                                                                                      \right.
\end{align*}
Therefore, we have $$\omega(f)=\begin{cases}
	   \left\lceil \frac{2kn}{k+1} \right\rceil, \textrm{if } n \equiv r \bmod (k+1) \textrm{ for some }  0<r<\frac{k+1}{2}; \\
	   \left\lceil \frac{2kn}{k+1} \right\rceil + 1, \mbox{otherwise}
\end{cases}$$
which implies
$$\gamma_{k-\text{SRD}}(P_2 \times P_n) \leqslant \begin{cases}
	   \left\lceil \frac{2kn}{k+1} \right\rceil, \textrm{if } n \equiv r \bmod (k+1) \textrm{ for some }  0<r<\frac{k+1}{2}; \\
	   \left\lceil \frac{2kn}{k+1} \right\rceil + 1, \mbox{otherwise}
\end{cases}$$

\underline{Lower bound}. Let $g:V(P_2\times P_n)\rightarrow S$ be an arbitrary $k$-SRD function defined on vertices of grid graph $P_2\times P_n$, $n\geqslant\frac{k}{2}$.  This function generates the decomposition $\{V_s:s\in S\}$ of $V(P_2\times P_n)$, such that $V_s:=\{a_{ij}\in V:g(a_{ij})=s\}$, $s\in S$. As we need to prove inequality
$$\omega(g)\geqslant \begin{cases}
	   \left\lceil \frac{2kn}{k+1} \right\rceil, \textrm{if } n \equiv r \bmod (k+1) \textrm{ for some }  0<r<\frac{k+1}{2}; \\
	   \left\lceil \frac{2kn}{k+1} \right\rceil + 1, \mbox{otherwise}. \end{cases}$$

To do so, it is sufficient to describe the appropriate discharging procedure in which weights of vertices from $V_2\cup V_3\cup V_4$ are used in order to charge weights of vertices from $V_0$ and accomplishing that, in the end, the following conditions are fulfilled:

\begin{itemize}
  \item Every vertex from  $V(P_2\times P_n)$ has a weight of at least $\frac{k}{k+1}$.
  \item In the case $n \equiv r \bmod (k+1)$, with $r=0$ or $r\geqslant\frac{k+1}{2}$, some vertices have ``excess'' weights such that the total sum $\Sigma_e$ of those
   excess weights is at least $1$.
\end{itemize}

Firstly, we introduce set $W:=\{a_{ij}\in V_3\cup V_4:N(a_{ij})\cap V_0\ne\emptyset\}$. If $W \neq \emptyset$ , then for every vertex $a_{ij}\in W$, we perform the following ``primitive'' discharging procedure:
\begin{enumerate}
  \item Vertices from $N(a_{ij})\cap V_0$ which are solely defended by the vertex $a_{ij}$ receive weight of $1$ from this vertex. Note that, in the case $a_{ij}\in V_3$ and $k>2$, there are at most $2$ vertices from $N(a_{ij})\cap V_0$ that can be charged this way.
  \item Some or all of the remaining vertices from $N(a_{ij})\cap V_0$ receive weight of $1$ from the remaining weight of vertex $a_{ij}$. This discharging procedure
  is terminated when all vertices from $N(a_{ij})\cap V_0$ are charged with weight $1$ or until the remaining weight of vertex $a_{ij}$ is $2$.
\end{enumerate}

Applying the previous discharging procedure, we obtain a $k-$SRD function $h:V(P_2\times P_n)\rightarrow S$ such that $\omega(h)=\omega(g)$. Moreover, this function provides labeling such
that neighbors of every vertex labeled with zero are using labels from $\{0,1,2\}$ i.e. $h=(V_0,V_1,V_2)$.

In the next step, we must ensure that every remaining zero is charged with a total weight of $\frac{k}{k+1}$ received from one or more of its neighboring twos. Let $\sigma$
be the total number of vertices labeled by zero left after the primitive discharging is executed. We consider the following cases:

\textbf{Case 1.} $\sigma\leqslant k$;

In this case, for each zero we can pick a unique two as its defender. This $1$--$1$ correspondence ensures that every such two can fully charge the corresponding zero with
the weight $\frac{k}{k+1}$. The remaining weight of each defender is equal to $2-\frac{k}{k+1}=\frac{k}{k+1}+\frac{2}{k+1}$, so each one of them has an excess weight of $\frac{2}{k+1}$.
Since the weight of any other non-zero labelled vertex is at least $1=\frac{k}{k+1}+\frac{1}{k+1}$, we have
$ \Sigma_e\geqslant\frac{2\sigma}{k+1}+\frac{2n-2\sigma}{k+1}=\frac{2n}{k+1}.$
Obviously, we have $\Sigma_e\geqslant 1$ for $n\geqslant\frac{k+1}{2}$, and $\Sigma_e<1$ otherwise.

\textbf{Case 2.} $\sigma> k$;

\begin{figure}[h!]
\adjustbox{scale=0.75,center}{
    	\begin{tikzcd}
    	 &	0  & 2 \arrow[l,dashrightarrow, bend right=50,"{\frac{k-1}{k+1}}"] \arrow[d,dashrightarrow, bend right=50,"{\frac{1}{k+1}}"]\arrow[r,dashrightarrow, swap, bend left=50,"{\frac{2}{k+1}}"]& 0 & 2 \arrow[l,dashrightarrow, bend right=50,"{\frac{k-3}{k+1}}"] \arrow[d,dashrightarrow, bend right=50,"{\frac{1}{k+1}}"] \arrow[r,dashrightarrow, swap, bend left=50,"{\frac{4}{k+1}}"] & \dotsc \arrow[r,dashrightarrow, swap, bend left=50,"{\frac{k-4}{k+1}}"]& 0 & 2 \arrow[r,dashrightarrow, swap, bend left=50,"{\frac{k-2}{k+1}}"] \arrow[d,dashrightarrow, bend right=50,"{\frac{1}{k+1}}"] \arrow[l,dashrightarrow, bend right=50,"{\frac{3}{k+1}}"]& 0 & 2 \arrow[r,dashrightarrow, swap, bend left=50,"{\frac{k}{k+1}}"] \arrow[d,dashrightarrow, bend right=50,"{\frac{1}{k+1}}"] \arrow[l,dashrightarrow, bend right=50,"{\frac{1}{k+1}}"]& 0\\
       0   &   2 \arrow[l,dashrightarrow, bend left=50,"{\frac{k}{k+1}}"] \arrow[u,dashrightarrow, swap, bend left=50,"{\frac{1}{k+1}}"] \arrow[r,dashrightarrow, swap, bend right=50,"{\frac{1}{k+1}}"] & 0  &   2 \arrow[l,dashrightarrow, bend left=50,"{\frac{k-2}{k+1}}"] \arrow[u,dashrightarrow, swap, bend left=50,"{\frac{1}{k+1}}"]  \arrow[r,dashrightarrow, swap, bend right=50,"{\frac{3}{k+1}}"] &  0 & \dotsc \arrow[l,dashrightarrow, bend left=50,"{\frac{k-4}{k+1}}"] & 2 \arrow[r,dashrightarrow, swap, bend right=50,"{\frac{k-3}{k+1}}"] \arrow[u,dashrightarrow, swap, bend left=50,"{\frac{1}{k+1}}"] \arrow[l,dashrightarrow, bend left=50,"{\frac{4}{k+1}}"]& 0 & 2 \arrow[r,dashrightarrow, swap, bend right=50,"{\frac{k-1}{k+1}}"] \arrow[u,dashrightarrow, swap, bend left=50,"{\frac{1}{k+1}}"] \arrow[l,dashrightarrow, bend left=50,"{\frac{2}{k+1}}"]& 0
    	\end{tikzcd}
    }
\caption{The discharging rule for even $k$}\label{fig:discharging_even_k_grid_2_n}

\end{figure}
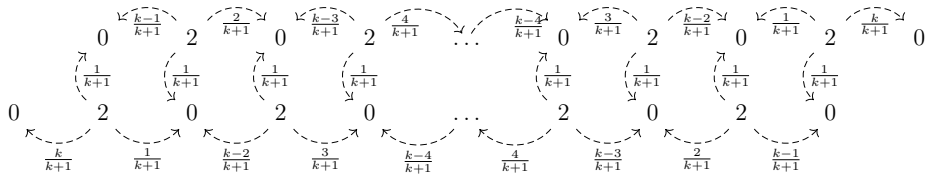

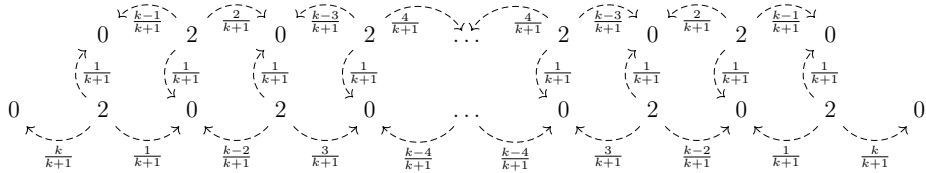
\begin{figure}[h!]
\adjustbox{scale=0.75,center}{
    	\begin{tikzcd}
    	 &	0  & 2 \arrow[l,dashrightarrow, bend right=50,"{\frac{k-1}{k+1}}"] \arrow[d,dashrightarrow, bend right=50,"{\frac{1}{k+1}}"]\arrow[r,dashrightarrow, swap, bend left=50,"{\frac{2}{k+1}}"]& 0 & 2 \arrow[l,dashrightarrow, bend right=50,"{\frac{k-3}{k+1}}"] \arrow[d,dashrightarrow, bend right=50,"{\frac{1}{k+1}}"] \arrow[r,dashrightarrow, swap, bend left=50,"{\frac{4}{k+1}}"] & \dotsc & 2 \arrow[l,dashrightarrow, bend right=50,"{\frac{4}{k+1}}"] \arrow[r,dashrightarrow, swap, bend left=50,"{\frac{k-3}{k+1}}"] \arrow[d,dashrightarrow, bend right=50,"{\frac{1}{k+1}}"]& 0  & 2 \arrow[r,dashrightarrow, swap, bend left=50,"{\frac{k-1}{k+1}}"] \arrow[d,dashrightarrow, bend right=50,"{\frac{1}{k+1}}"] \arrow[l,dashrightarrow, bend right=50,"{\frac{2}{k+1}}"]& 0 \\
       0   &   2 \arrow[l,dashrightarrow, bend left=50,"{\frac{k}{k+1}}"] \arrow[u,dashrightarrow, swap, bend left=50,"{\frac{1}{k+1}}"] \arrow[r,dashrightarrow, swap, bend right=50,"{\frac{1}{k+1}}"] & 0  &   2 \arrow[l,dashrightarrow, bend left=50,"{\frac{k-2}{k+1}}"] \arrow[u,dashrightarrow, swap, bend left=50,"{\frac{1}{k+1}}"]  \arrow[r,dashrightarrow, swap, bend right=50,"{\frac{3}{k+1}}"] &  0 & \dotsc \arrow[l,dashrightarrow, bend left=50,"{\frac{k-4}{k+1}}"] \arrow[r,dashrightarrow, swap, bend right=50,"{\frac{k-4}{k+1}}"]& 0  & 2 \arrow[l,dashrightarrow, bend left=50,"{\frac{3}{k+1}}"] \arrow[r,dashrightarrow, swap, bend right=50,"{\frac{k-2}{k+1}}"] \arrow[u,dashrightarrow, swap, bend left=50,"{\frac{1}{k+1}}"]& 0  & 2 \arrow[r,dashrightarrow, swap, bend right=50,"{\frac{k}{k+1}}"] \arrow[u,dashrightarrow, swap, bend left=50,"{\frac{1}{k+1}}"] \arrow[l,dashrightarrow, bend left=50,"{\frac{1}{k+1}}"]& 0
    	\end{tikzcd}
    }
\caption{The discharging rule for odd $k$} \label{fig:discharging_odd_k_grid_2_n}
\end{figure}

If any $k+1$ zeros cannot be covered with less of $k+1$ twos, then we use a similar reasoning as in Case 1. Let us denote this type of pattern with $P$. In type pattern $P$, every defender remains with weight $\frac{2}{k+1}$, implying that the total excess weight of this $2k+2$ vertices is equal to $2$. Besides type pattern $P$, there can exist a configuration in which exactly $k+1$ zeros are covered with exactly $k$ twos (type pattern $P'$), and a configuration in which exactly $k+2$ zeros are covered with exactly $k$ twos (type pattern $P''$). Concerning the discharging procedure, the most challenging configuration is the type pattern $P''$. This pattern consists of a zig-zag path obtained by alternating $k+2$ zeros and $k$ twos. In this situation, the memory-based discharging procedure within this path is described in Figures~\ref{fig:discharging_even_k_grid_2_n}-- \ref{fig:discharging_odd_k_grid_2_n}.

More precisely, $m-$th leftmost two in the observed path charge its left zero neighbor with weight $\frac{k-m+1}{k+1}$, charge its right zero neighbor with weight
$\frac{m}{k+1}$, and charge its up or down zero neighbor with weight $\frac{1}{k+1}$. As a result of applying the described procedure, each vertex within the considered path has a weight of exactly $\frac{k}{k+1}$, implying that this path does not contribute to the increase of total excess weight $\Sigma_e$.

One example of type pattern $P'$ can be derived from the described type pattern $P''$ by disregarding one of the outermost zeros. The described discharging procedure is also valid for this type of pattern with the difference that, in this case, type pattern $P'$ contributes to the increase of total excess weight $\Sigma_e$ by $\frac{k}{k+1}$. A discharging procedure for other configurations of type pattern $P'$ (exactly $k+1$ zeros covered by exactly $k$ twos) can be constructed in a similar manner. It can be shown that any pattern of this
type contributes to the increase of total excess weight $\Sigma_e$ by at least $\frac{k}{k+1}$.

Finally, we can use previous observations to characterize the condition $\Sigma_e<1$. Let $n=(k+1)l+r$, with $l\geqslant 1$ and $0\leqslant r\leqslant k$. The existence of
at least one type pattern $P$ or the existence of at least two disjoint type patterns $P'$ would imply $\Sigma_e\geqslant 1$. The same conclusion holds when one copy of type pattern $P'$ is presented, since on the other side of the corresponding pattern $P''$ exists a vertex with the excess weight of at least  $\frac{1}{k+1}$ (note that
this situation occurs in the case $r=0$). So, the only remaining possibility is to concatenate several copies of pattern $P''$, including their upside-down reflection counterparts. Since $2n=(2k+2)l+2r$, $0<r\leqslant 2k$, we can use at most $l$ copies of
pattern $P''$, so the residue of the grid graph (on both sides) consists of $2r$ vertices. If $\sigma'$ is the total number of zeros in this residue, then $\sigma'\leqslant k$
(since $2r\leqslant 2k$), so we have at least $\sigma'$ twos in this residue which are able to defend those zeros. If residue contains a non-defender vertex from $V_2$, then this vertex all by itself has the excess weight greater than $1$. Otherwise, after performing the discharging procedure of this residue, we get
$\Sigma_e=\frac{2\sigma'}{k+1}+\frac{2r-2\sigma'}{k+1}=\frac{2r}{k+1},$  and condition $\Sigma_e<1$ is equivalent to condition $0<r<\frac{k+1}{2}$.

Taking all these facts together, we get the inequality
$$ \omega(g)=\omega(h)\geqslant \begin{cases}
	    \frac{2kn}{k+1} , \textrm{if } n \equiv r \bmod (k+1) \textrm{ for some }  0<r<\frac{k+1}{2}; \\
	    \frac{2kn}{k+1} + 1, \mbox{otherwise.}
\end{cases}  $$

Since $g$ was an arbitrary $k-SRD$ function on $P_2\times P_n$, we conclude
$$ \gamma_{k-\text{SRD}}(P_2 \times P_n) \geqslant \begin{cases}
	   \lceil \frac{2kn}{k+1} \rceil, \textrm{if } n \equiv r \bmod (k+1) \textrm{ for some }  0<r<\frac{k+1}{2}; \\
	   \lceil \frac{2kn}{k+1} \rceil + 1, \mbox{otherwise.}
\end{cases}  $$
\qed

  \section{Conclusions}

This paper addresses the $k$-strong Roman domination problem, an extension of the widely studied Roman domination problem on graphs, which traditionally involves single attacks on each vertex. Unlike the original problem, this generalized version considers multiple simultaneous attacks, that occur on any $k$-combination of graph vertices. This problem provides a more realistic depiction of defense strategies on real networks.

This work is pioneering in its theoretical exploration of the $k$-strong Roman domination problem. Notably, it establishes the exact bound for the $k$-SRD problem across several fundamental graph classes, including complete graphs, paths, cycles, wheels, complete bipartite graphs, and grid graphs $2 \times n$, in the case of arbitrary $k\geqslant 2$.   For fixed $k=2$, general cubic graphs are considered, providing a tight lower bound is by means of the discharging approach.  Moreover, two cubic graph classes are identified, belonging to the class of convex polytopes, whose exact bounds attain the general lower bound for cubic graphs. Finally, for the convex polytope $A_n$, an upper bound on the $k$-SRD number is established.

For future work, one could prove the exact bound on $A_n$  which remains here unresolved and for which we state a hypothesis that it is equal to the provided upper bound. Additionally, we believe that a tight upper bound on cubic graphs is not far from the lower bound we provided in this work -- in particular, we state a hypothesis that it is exactly equal to  $\left \lceil \frac{2}{3} n\right \rceil +2$. Also, identifying the exact bound in the case of general grid graphs $m \times n$ for arbitrary $k\geqslant 2$ is another open question.    Detecting interesting relations between the $k$-strong Roman domination and the strong Roman domination problem, which the both of them consider simultaneous attacks, is worth an investigation in the near future.

\section*{Competing interests}
The authors declare no competing interests.

\section*{Acknowledgments}
All authors  are supported by the bilateral project between Austria and Bosnia and Herzegovina funded by the Ministry of Civil Affairs of Bosnia
and Herzegovina under grant no. 1259074.

   \bibliographystyle{abbrv}	
 \bibliography{bib}

\begin{thebibliography}{10}

\bibitem{alvarez2017strong}
M.~{\'A}lvarez-Ruiz, T.~Mediavilla-Gradolph, S.~M. Sheikholeslami,
  J.~Valenzuela-Tripodoro, and I.~G. Yero.
\newblock On the strong {R}oman domination number of graphs.
\newblock {\em Discrete Applied Mathematics}, 231:44--59, 2017.

\bibitem{baca1988labelings}
M.~Ba{\v{c}}a.
\newblock Labelings of 2 classes of convex polytopes.
\newblock {\em Utilitas Mathematica}, 34:24--31, 1988.

\bibitem{balasundaram2006graph}
B.~Balasundaram and S.~Butenko.
\newblock Graph domination, coloring and cliques in telecommunications.
\newblock {\em Handbook of optimization in telecommunications}, pages 865--890,
  2006.

\bibitem{chang1998algorithmic}
G.~J. Chang.
\newblock Algorithmic aspects of domination in graphs.
\newblock {\em Handbook of Combinatorial Optimization: Volume1--3}, pages
  1811--1877, 1998.

\bibitem{chellali2020roman}
M.~Chellali, N.~Jafari~Rad, S.~M. Sheikholeslami, and L.~Volkmann.
\newblock Roman domination in graphs.
\newblock {\em Topics in domination in graphs}, pages 365--409, 2020.

\bibitem{chellali2021varieties}
M.~Chellali, N.~J. Rad, S.~Sheikholeslami, and L.~Volkmann.
\newblock Varieties of {R}oman domination.
\newblock {\em Structures of domination in graphs}, pages 273--307, 2021.

\bibitem{church1974maximal}
R.~Church and C.~ReVelle.
\newblock The maximal covering location problem.
\newblock In {\em Papers of the regional science association}, volume~32, pages
  101--118. Springer-Verlag Berlin/Heidelberg, 1974.

\bibitem{cockayne1978domination}
E.~Cockayne.
\newblock Domination of undirected graphs—a survey.
\newblock In {\em Theory and Applications of Graphs: Proceedings, Michigan May
  11--15, 1976}, pages 141--147. Springer, 1978.

\bibitem{cockayne2004roman}
E.~J. Cockayne, P.~A. Dreyer~Jr, S.~M. Hedetniemi, and S.~T. Hedetniemi.
\newblock Roman domination in graphs.
\newblock {\em Discrete mathematics}, 278(1-3):11--22, 2004.

\bibitem{filipovic2022solving}
V.~Filipovi{\'c}, D.~Mati{\'c}, and A.~Kartelj.
\newblock Solving the signed {R}oman domination and signed total {R}oman
  domination problems with exact and heuristic methods.
\newblock {\em arXiv preprint arXiv:2201.00394}, 2022.

\bibitem{goddard2013independent}
W.~Goddard and M.~A. Henning.
\newblock Independent domination in graphs: A survey and recent results.
\newblock {\em Discrete Mathematics}, 313(7):839--854, 2013.

\bibitem{gupta2013domination}
P.~Gupta.
\newblock Domination in graph with application.
\newblock {\em Indian J. Res}, 2(3):115--117, 2013.

\bibitem{henning2009survey}
M.~A. Henning.
\newblock A survey of selected recent results on total domination in graphs.
\newblock {\em Discrete Mathematics}, 309(1):32--63, 2009.

\bibitem{henning2003defending}
M.~A. Henning and S.~T. Hedetniemi.
\newblock Defending the {R}oman empire—a new strategy.
\newblock {\em Discrete Mathematics}, 266(1-3):239--251, 2003.

\bibitem{kammerling2009roman}
K.~Kammerling and L.~Volkmann.
\newblock Roman k-domination in graphs.
\newblock {\em Journal of the Korean Mathematical Society}, 46(6):1309--1318,
  2009.

\bibitem{kartelj2021roman}
A.~Kartelj, M.~Grbi{\'c}, D.~Mati{\'c}, and V.~Filipovi{\'c}.
\newblock The {R}oman domination number of some special classes of
  graphs-convex polytopes.
\newblock {\em Applicable Analysis and Discrete Mathematics}, 15(2):393--412,
  2021.

\bibitem{liu2013roman}
C.-H. Liu and G.~J. Chang.
\newblock Roman domination on strongly chordal graphs.
\newblock {\em Journal of Combinatorial Optimization}, 26(3):608--619, 2013.

\bibitem{liu2020k}
Z.~Liu, X.~Li, and A.~Khojandi.
\newblock On the k-strong {R}oman domination problem.
\newblock {\em Discrete Applied Mathematics}, 285:227--241, 2020.

\bibitem{lu2010survey}
G.~Lu, M.-T. Zhou, Y.~Tang, Z.-Q. Wu, G.-Y. Qiu, and L.~Yuan.
\newblock A survey on exact algorithms for dominating set related problems in
  arbitrary graphs.
\newblock {\em Chinese Journal of Computers}, 33(6):1073--1087, 2010.

\bibitem{milenkovic2011dominating}
T.~Milenkovi{\'c}, V.~Memi{\v{s}}evi{\'c}, A.~Bonato, and N.~Pr{\v{z}}ulj.
\newblock Dominating biological networks.
\newblock {\em PloS one}, 6(8):e23016, 2011.

\bibitem{macdougall2006vertex}
M.~Miller, M.~Ba{\v{c}}a, and J.~A. MacDougall.
\newblock Vertex-magic total labeling of generalized petersen graphs and convex
  polytopes.
\newblock {\em Journal of Combinatorial Mathematics and Combinatorial
  Computing}, 58:89--99, 2006.

\bibitem{pagourtzis2002server}
A.~Pagourtzis, P.~Penna, K.~Schlude, K.~Steinh{\"o}fel, D.~S. Taylor, and
  P.~Widmayer.
\newblock Server placements, {R}oman domination and other dominating set
  variants.
\newblock In {\em Foundations of Information Technology in the Era of Network
  and Mobile Computing: IFIP 17th World Computer Congress—TC1 Stream/2nd IFIP
  International Conference on Theoretical Computer Science (TCS 2002) August
  25--30, 2002, Montr{\'e}al, Qu{\'e}bec, Canada}, pages 280--291. Springer,
  2002.

\bibitem{rao2021survey}
Y.~Rao, R.~Chen, P.~Wu, H.~Jiang, and S.~Kosari.
\newblock A survey on domination in vague graphs with application in
  transferring cancer patients between countries.
\newblock {\em Mathematics}, 9(11):1258, 2021.

\bibitem{schoning1986complete}
U.~Sch{\"o}ning.
\newblock Complete sets and closeness to complexity classes.
\newblock {\em Mathematical Systems Theory}, 19(1):29--41, 1986.

\bibitem{shao2018discharging}
Z.~Shao, P.~Wu, H.~Jiang, Z.~Li, J.~{\v{Z}}erovnik, and X.~Zhang.
\newblock Discharging approach for double {R}oman domination in graphs.
\newblock {\em IEEE Access}, 6:63345--63351, 2018.

\end{thebibliography}

\end{document}